\documentclass{aims}

\usepackage{txfonts}
\def\typeofarticle{Type of article}
\def\currentvolume{8}
\def\currentissue{x}
\def\currentyear{2023}

\def\ppages{xxx--xxx}
\def\DOI{10.3934/math.2023xxx}
\def\Received{January 2023}
\def\Revised{February 2023}
\def\Accepted{March 2023}
\def\Published{March 2023}

\numberwithin{equation}{section}

\makeatletter
\renewcommand{\@biblabel}[1]{#1\hfill \hspace{-0.2cm}}
\makeatother

\usepackage{pifont}
\usepackage{color}
\usepackage{amssymb}
\usepackage{mathrsfs}
\usepackage{amsmath, amsfonts, enumerate}
\usepackage{graphicx}
\usepackage{subcaption}
\usepackage{verbatim}
\usepackage{bm}
\usepackage{euscript}
\usepackage{pgflibraryarrows}
\usepackage{pgflibrarysnakes}
\usepackage{subeqnarray}
\usepackage{cases}

\renewcommand{\l}{\lambda}
\newcommand{\dd}{\displaystyle}

\renewcommand{\l}{\lambda}
\newcommand{\R}{\mathbb{R}}

\newcommand{\s}{\sigma}
\newcommand{\x}{\xi}
\newcommand{\am}{\operatorname{am}}
\newcommand{\cn}{\operatorname{cn}}
\newcommand{\sn}{\operatorname{sn}}
\newcommand{\dn}{\operatorname{dn}}

\newcommand{\cnap}[1]{\cn^{#1}(\l \x , m)}

\begin{document}

\title{Exact Jacobi elliptic solutions of some models for the interaction of long and short waves.}

\author{%
Bruce Brewer\affil{1},
  Jake Daniels\affil{1}
  and
  Nghiem V. Nguyen \affil{1,}\corrauth
}

\shortauthors{Bruce Brewer, Jake Daniels and Nghiem V. Nguyen}

\address{%
 \addr{\affilnum{1}}{Department of Mathematics and Statistics, Utah State University, 3900 Old Main Hill, Logan, UT 84322, USA}
}

\corraddr{Email: nghiem.nguyen@usu.edu; Tel: +1-435-797-2819.
}

\begin{abstract}
Some systems were recently put forth by Nguyen \textit{et. al.} as models for studying the interaction of long and short waves in dispersive media.  These systems were shown to possess synchronized Jacobi elliptic solutions as well as synchronized solitary wave solutions under certain constraints, \textit{i.e.}, vector solutions where the two components are proportional to one another.    In this paper, the exact periodic traveling wave solutions to these systems in general are found to be given by Jacobi elliptic functions.  Moreover, these cnoidal wave solutions are unique.  Thus, the explicit synchronized solutions under some conditions obtained by Nguyen \textit{et. al.} are also indeed unique.
\end{abstract}

\keywords{Periodic solutions, Cnoidal solutions, NLS-equation, KdV-equation, BBM-equation, NLS-KdV system.
\newline
\textbf{Mathematics Subject Classification:} 35A16, 35A24, 35B10. }

\maketitle

\section{Introduction}
The following four systems, termed Schr\"odinger KdV-KdV, Schr\"odinger BBM-BBM, Schr\"odinger KdV-BBM and Schr\"odinger BBM-KdV, respectively,
\begin{equation}\label{system1}
\left\{
\begin{matrix}
\begin{split}
&\frac{\partial u}{\partial t}+\mu_0 \frac{\partial u}{\partial x}+a_0\frac{\partial^3 u}{\partial x^3} + i  b \frac{\partial^2 u}{\partial x^2}=-\frac{\partial (uv)}{\partial x}-i\mu_1 uv,\\
&\frac{\partial v}{\partial t}+\frac{\partial v}{\partial x}
  +v\frac{\partial v}{\partial x}+c\frac{\partial^3 v}{\partial x^3}= -\frac{1}{2}\frac{\partial |u|^2}{\partial x};
\end{split}
\end{matrix}
\right.
\end{equation}

\begin{equation}\label{system2}
\left\{
\begin{matrix}
\begin{split}
&\frac{\partial u}{\partial t}+\mu_0 \frac{\partial u}{\partial x}-a_1\frac{\partial^3 u}{\partial x^2 \partial t} + i  b \frac{\partial^2 u}{\partial x^2}=-\frac{\partial (uv)}{\partial x}-i\mu_1 uv,\\
&\frac{\partial v}{\partial t}+\frac{\partial v}{\partial x}
  +v\frac{\partial v}{\partial x}-c\frac{\partial^3 v}{\partial x^2\partial t}= -\frac{1}{2}\frac{\partial |u|^2}{\partial x};
\end{split}
\end{matrix}
\right.
\end{equation}
\begin{equation}\label{system3}
\left\{
\begin{matrix}
\begin{split}
&\frac{\partial u}{\partial t}+\mu_0 \frac{\partial u}{\partial x}+a_0\frac{\partial^3 u}{\partial x^3} + i  b \frac{\partial^2 u}{\partial x^2}=-\frac{\partial (uv)}{\partial x}-i\mu_1 uv,\\
&\frac{\partial v}{\partial t}+\frac{\partial v}{\partial x}
  +v\frac{\partial v}{\partial x}-c\frac{\partial^3 v}{\partial x^2\partial t}= -\frac{1}{2}\frac{\partial |u|^2}{\partial x};
\end{split}
\end{matrix}
\right.
\end{equation}
and
\begin{equation}\label{system4}
\left\{
\begin{matrix}
\begin{split}
&\frac{\partial u}{\partial t}+\mu_0 \frac{\partial u}{\partial x}-a_1\frac{\partial^3 u}{\partial x^2\partial t} + i  b \frac{\partial^2 u}{\partial x^2}=-\frac{\partial (uv)}{\partial x}-i\mu_1 uv,\\
&\frac{\partial v}{\partial t}+\frac{\partial v}{\partial x}
  +v\frac{\partial v}{\partial x}+c\frac{\partial^3 v}{\partial x^3}= -\frac{1}{2}\frac{\partial |u|^2}{\partial x};
\end{split}
\end{matrix}
\right.
\end{equation}
were recently advocated in \cite{LN,LN2} (see also \cite{DNS}) as more suitable models for studying the interaction of long and short waves in dispersive media due to their consistent derivation when compared to the nonlinear Schr\"odinger-Korteweg-de Vries system \cite{KSK}
\begin{equation}
\left\{
\begin{matrix}
\begin{split}
i u_t + u_{xx} + a|u|^2u &= -buv,\\
v_t + c vv_x + v_{xxx} &= -\frac{b}{2}(|u|^2)_x.
\end{split}
\end{matrix}
\right.
\label{NLS-KdV}
\end{equation}  Here, the function $u(x,t)$ is a complex-valued function while $v(x,t)$ is a real-valued function, and  $x,t \in \mathbb R$ where $\mu_0,\mu_1, a_0,a_1,b$ and $c$ are real constants with $\mu_0,\mu_1, a_0, a_1,c>0$. 
 For detailed discussion on these systems, we refer our readers to the papers \cite{LN,LN2,DNS}.  

A traveling-wave solution to the above four systems is a vector solution $\big(u(x,t),v(x,t)\big)$ of the form 
\begin{equation}
\label{travelingwavesolution}
u(x,t) = e^{i\omega t} e^{iB(x-\sigma t)}f(x-\sigma t),\ \ \ \ \ \ \ v(x,t) = g(x-\sigma t),
\end{equation} where $f$ and $g$ are smooth, real-valued functions with speed $\sigma>0$, and phase shifts $B,\omega \in \mathbb R$.  Substituting the traveling-wave ansatz \eqref{travelingwavesolution} into the four systems and separating the real and imaginary parts, the following associated ODE systems are obtained
\begin{equation}\label{associatedODE1}
\left\{
\begin{matrix}
\begin{split}
& f' g + f g' + a_0 f''' + (\mu_0 - \sigma - 3 a_0 B^2 - 2 b B) f' = 0, \\
&    (B + \mu_1) f g + (3 a_0 B + b) f'' + (\omega + B \mu_0 - B \sigma - a_0 B^3 - b B^2) f = 0, \\
 &   f f' + g g' + c g'''+ (1 - \sigma) g' = 0;
\end{split}
\end{matrix}
\right.
\end{equation}
\begin{equation}\label{associatedODE2}
\left\{
\begin{matrix}
\begin{split}
 &  f' g + f g' + a_1 \sigma f''' + (\mu_0 + 2 a_1 B \omega - 3 a_1 B^2 \sigma - \sigma - 2 b B) f' = 0, \\
   & (B + \mu_1) f g + (3 a_1 B \sigma + b - a_1 \omega) f'' + (\omega + B \mu_0 + a_1 B^2 \omega - a_1 B^3 \sigma - B \sigma - b B^2) f = 0, \\
  &  ff' + g g' + c \sigma g''' + (1 - \sigma) g' = 0;
\end{split}
\end{matrix}
\right.
\end{equation}
\begin{equation}\label{associatedODE3}
\left\{
\begin{matrix}
\begin{split}
&  f' g + f g' + a_0 f''' + (\mu_0 - \sigma - 3 a_0 B^2 - 2 b B) f' = 0, \\
  &  (B + \mu_1) f g + (3 a_0 B + b) f'' + (\omega + B \mu_0 - B \sigma - a_0 B^3 - b B^2) f = 0, \\
   & f f' + g g' + c \sigma g''' + (1 - \sigma) g' = 0;
\end{split}
\end{matrix}
\right.
\end{equation}
and
\begin{equation}\label{associatedODE4}
\left\{
\begin{matrix}
\begin{split}
&  f' g + f g' + a_1 \sigma f''' + (\mu_0 + 2 a_1 B \omega - 3 a_1 B^2 \sigma - \sigma - 2 b B) f' = 0, \\
 &   (B + \mu_1) f g + (3 a_1 B \sigma + b - a_1 \omega) f'' + (\omega + B \mu_0 + a_1 B^2 \omega - a_1 B^3 \sigma - B \sigma - b B^2) f = 0, \\
  &  f f' + g g' + c g'''+ (1 - \sigma) g' = 0.
\end{split}
\end{matrix}
\right.
\end{equation}
We refer to semi-trivial solutions as solutions where at least one component is a constant (possibly zero).   Of course, the trivial solution $(0,0)$ is always a solution. In the case when $f$ is a constant multiple of $g$, the vector solution is termed synchronized solution.  Among the traveling-wave solutions, attention is often given to the solitary-wave and  periodic solutions due to the roles they sometimes play in the evolution equations.  Solitary waves are smooth traveling-wave solutions that are symmetric around a \textit{single} maximum and rapidly decay to zero away from the maximum while periodic solutions are self-explanatory. Even though less common, the term solitary waves are also sometimes used to describe traveling-wave solutions that are symmetric around a single maximum but that approach nonzero constants as $\xi\rightarrow \pm \infty$. 

The topic of existence of synchronized traveling-wave solutions to these four systems has been addressed previously (\cite{NLB}).  Notice that when $f$ is a constant multiple of $g$, \textit{i.e.} $u=Av$ for some proportional constant $A$, the three equations in each of the four associated ODEs \eqref{associatedODE1}-\eqref{associatedODE4} can be collapsed into four single equations of the form
\begin{equation}\label{genericODE}
f'^2 = k_3f^3 +k_2f^2 + k_1 f + k_0
\end{equation} under certain constraints.
  In \cite{NLB}, it was shown that the systems possess synchronized solitary waves with the usual hyperbolic $\operatorname{sech^2}$-profile typical of dispersive equations.  In \cite{BLN},  a novel approach was first employed to establish the existence of periodic traveling-wave solutions for these systems, namely the topological degree theory for positive operators that was introduced by Krasnosel'skii \cite{K1,K2} and used in several different models \cite{C,CCN,N}.  The explicit synchronized periodic solutions $u = Av$ where $v$ is given by the Jacobi elliptic function
\begin{equation}\label{Jacobiform}
v(x-\sigma t) :\equiv v(\xi) = C_0 + C_2  \cn^2(\alpha \xi + \beta,m)
\end{equation}
were then obtained by demanding the coefficients in each of the four cases to satisfy certain constraints.  (A brief description of the Jacobi elliptic functions is recalled below.)   Neither approach, however, guarantees uniqueness of the periodic solutions obtained due to several factors such as the form of $v$ as a priori assumption because of \eqref{genericODE} as well as the nature of the topological degree theory approach. 

It is worth to point out that explicit solitary wave solutions have been found for another system \cite{C1,C2}, the $abcd$-system
\begin{equation}\label{abcd-system}
    \begin{cases}
        \eta_t + w_x + (w\eta)_x + aw_{xxx} - b\eta_{xxt} = 0 \,,\\
        w_t + \eta_x +ww_x + c\eta_{xxx} - dw_{xxt} = 0 \,,
    \end{cases}
\end{equation}
where $a, b, c,$ and $d$ are real constants satisfying
\begin{align*}
    a + b &= \frac{1}{2}(\theta^2 - \frac{1}{3}), \\
    c + d &= \frac{1}{2}(1 - \theta^2) \geq 0, \\
    a + b + c + d &= \frac{1}{3},
\end{align*}
and $\theta \in [0,1]$. This system is used to model small-amplitude,
long wavelength, gravity waves on the surface of water 
(\cite{BCS1,BCS2}). Here $\eta(x,t)$ and $w(x,t)$ are real valued
functions and $x, t \in \R$.  However, the existence of periodic traveling wave solutions for this system is still not well understood.  The only result that we are aware of is for the special case when $a=c=0$ and $b=d=1/6$, where the solutions are given in term of the Jacobi elliptic cnoidal function \cite{CCN}.

The manuscript is organized as follows.  In Section 2, some facts about the Jacobi elliptic functions are reviewed and the results are summarized.
 In Section 3, the explicit cnoidal solutions to the four systems are established, and how these solutions limit to the solitary-wave solutions is analysed.  Section 5 is devoted to discussion of the obtained results.  To preserve the self-completeness without affecting the flow of the paper, some tedious formulae and expressions are delegated to the Appendix.

\section{Preliminaries and statement of results}
For the readers' convinience, some notions of the Jacobi elliptic functions are briefly recalled here.
Let
\begin{equation*}
v=\int_0^{\phi} \frac{1}{\sqrt{1-m^2 \sin^2 t}}dt, \ \ \ \ \ \mbox{for $0\leq m\leq 1$}.
\end{equation*}
Then $v =F(\phi, m)$ or equivalently, $\phi = F^{-1}(v,m)= \am(v,m)$ which is the Jacobi amplitude.  The two basic Jacobi elliptic functions $\cn(v,m)$ and $\sn(v,m)$ are defined as
\begin{equation*}
\sn(v,m) = \sin(\phi) = \sin\big(F^{-1}(v,m)\big) \ \ \ \ \ \mbox{and} \ \ \ \ \ \cn(v,m)= \cos(\phi) = \cos\big(F^{-1}(v,m)\big),
\end{equation*}
where $m$ is referred to as the Jacobi elliptic modulus. These functions are generalizations of the trigonometric and hyperbolic functions which satisfy
\begin{equation*}
\begin{matrix}
\sn(v,0) =\sin(v), \ \ \ \ \ \ \ \cn(v,0) =\cos(v),\\
\cn(v,1) = \operatorname{sech} (v), \ \ \ \ \ \sn(v,1) =\tanh (v).
\end{matrix}
\end{equation*}
We recall the following relations between these functions:
\begin{equation*}
\left\{
\begin{matrix}
\begin{split}
&\sn^2(\lambda \xi,m) =1-\cn^2 (\lambda \xi,m),\\
&\dn^2(\lambda \xi,m) =1 - m^2 + m^2 \cn^2(\lambda \xi,m),\\
&\frac{d}{d\xi}\cn(\lambda \xi,m)  = -\lambda \sn(\lambda \xi,m) \dn(\lambda \xi,m),\\
&\frac{d}{d\xi} \sn(\lambda \xi,m) = \lambda \cn(\lambda \xi,m) \dn(\lambda \xi,m),\\
&\frac{d}{d\xi} \dn(\lambda \xi,m) = -m^2 \lambda \cn(\lambda \xi,m) \sn(\lambda \xi,m).
\end{split}
\end{matrix}
\right.
\end{equation*}
In this manuscript, the existence of periodic traveling-wave solutions to the above four associated ODE systems \eqref{associatedODE1}-\eqref{associatedODE4} in general is analyzed.  The periodic traveling-wave solutions sought here are given by
\begin{equation}\label{eqn:FG}
    f(\x) = \sum_{r = 0}^n d_r \cnap{r} \quad \text{and} \quad g(\x) = \sum_{r = 0}^n h_r \cnap{r}\,,
\end{equation}
where $d_r,h_r \in \mathbb R$, $\lambda > 0$, and $0 \leq m \leq 1$.  Using the above relations, the following is revealed
\begin{equation}\label{derivatives}
\left\{
\begin{matrix}
\begin{split}
&\frac{d}{d\xi}\cn^r  = -r\lambda \cn^{r-1} \sn \dn,\\
&\frac{d^2}{d\xi^2} \cn^r = -r \lambda^2\big[(r+1) m^2 \cn^{r+2} + r (1-2m^2) \cn^r + (r-1) (m^2-1) \cn^{r-2}\big],\\
&\frac{d^3}{d\xi^3} \cn^r = r \lambda^3 \sn \dn \big[(r+1) (r+2) m^2 \cn^{r+1} + r^2 (1-2m^2) \cn^{r-1} + (r-1) (r-2) (m^2-1) \cn^{r-3}\big]
\end{split}
\end{matrix}
\right.
\end{equation}
where the argument $(\lambda \xi,m)$ has been dropped for clarity reason.  Notice that each of the above four associated ODE systems \eqref{associatedODE1}-\eqref{associatedODE4} involves three equations.  Plugging \eqref{derivatives} into these systems, the following generic form is obtained
\begin{equation}\label{genericsystem}
\left\{
\begin{matrix}
\begin{split}
&\sn(\lambda \xi,m) \dn(\lambda \xi,m) \sum_{q=0}^{2n-1} k_{1,q} \cn^q (\lambda \xi,m) =0,\\
& \sum_{q=0}^{2n} k_{2,q} \cn^q(\lambda \xi,m) =0,\\
&\sn(\lambda \xi,m) \dn(\lambda \xi,m) \sum_{q=0}^{2n-1} k_{3,q}\cn^q(\lambda \xi,m) =0
\end{split}
\end{matrix}
\right.
\end{equation}
where the subscripts $j$ and $q$ in the coefficient $k_{j,q}$ indicate the equation and the power on the cnoidal function $\cn$, respectively.  Notice that as \eqref{genericsystem} must hold true for all $(\lambda \xi,m)$, it must be the case that $k_{j,q}=0$ for each $j$ and $q$.  Moreover, from the third equation in all four systems, the sum $(ff'+gg')$ contributes the highest order term of $\cn^{2n-1}$. While the next highest order term is from $g'''$, which is $\cn^{n+1}$.  By balancing these highest order terms, it reveals that when $n\geq 3$, the highest order term is
$$k_{3,2n-1} \cn^{2n-1} = -n \lambda (d_n^2 + h_n^2) \cn^{2n-1}.$$
Since $\lambda, n>0$, requiring $k_{3,2n-1}=0$ implies that $d_n=h_n=0$, holding true for all $n\geq 3$.  Thus, the periodic traveling-wave ansatz \eqref{eqn:FG} reduces to
\begin{equation}\label{cnoidalsolutions}
f(\x)= d_0 +d_1 \cn(\lambda \xi,m) +d_2 \cn^2(\lambda \xi,m) \quad \text{and} \quad g(\x) = h_0 + h_1 \cn(\lambda \xi,m) + h_2 \cn^2(\lambda \xi,m).
\end{equation}

Next, by demanding all the coefficients $k_{j,q}=0$, a set of 13 equations is obtained for each of the four systems involving 11 unknowns  $d_i, h_i, B, \lambda, \omega, \sigma$ and $m$ with $i=0,1,2$.  (Equations \eqref{coeffKdVKdV}-\eqref{coeffBBMKdV}.)
For the Schr\"odinger KdV-KdV and Schr\"odinger BBM-BBM, the first and last equations in \eqref{associatedODE1} and \eqref{associatedODE2}, respectively, further yield $d_1 = h_1 = 0$.  In particular, the only non-trivial periodic solutions for the systems \eqref{system1} and \eqref{system2} are of the form
\begin{equation}\label{system1and2}
f(\x)= d_0 +d_2 \cn^2(\lambda \xi,m) \quad \text{and} \quad g(\x) = h_0 + h_2 \cn^2(\lambda \xi,m).
\end{equation} 
Under these conditions, the sets of 13 equations involving 11 unknowns, (equations \eqref{coeffKdVKdV} and \eqref{coeffBBMBBM}), reduce to sets of 7 equations with 9 unknowns.  
Similarly, for the Schr\"odinger KdV-BBM system and the Schr\"odinger BBM-KdV sytem, the first and last equations in \eqref{associatedODE3} and \eqref{associatedODE4}, respectively, reveal that $h_1 = 0$. Additionally, when substituting $h_1=0$ into \eqref{coeffKdVBBM} and \eqref{coeffBBMKdV}, the coefficients $k_{3,2}$ and $k_{3,0}$ in both systems, require that either $d_1=0$ or $d_0=d_2=0$.
When $d_1 = h_1 = 0,$ we have solutions of the form \eqref{system1and2}, where the sets of 13 equations involving 11 unkowns, (equations \eqref{coeffKdVBBM} and \eqref{coeffBBMKdV}), reduce to 7 equations with 9 unknowns.
When $d_0 = d_2 = h_1 = 0,$ we have that the only non-trivial periodic solutions for the systems \eqref{system3} and \eqref{system4} are of the form
\begin{equation*}
f(\x)= d_1 \cn(\lambda \xi,m) \quad \text{and} \quad g(\x) = h_0 +  h_2 \cn^2(\lambda \xi,m)
\end{equation*} 
in which case the sets of 13 equations involving 11 unknowns, (equations \eqref{coeffKdVBBM} and \eqref{coeffBBMKdV}), reduce to sets of 6 equations with 8 unknowns. 

The exact, explicit periodic traveling-wave solutions to  the four systems \eqref{system1}-\eqref{system4} could then be established by solving those reduced nonlinear systems with the help of the software Maple.  As there are two degrees of freedom, in principle any pair of two unknowns can be chosen as ``free parameters" so long as solutions can be found consistently.  In most physical situations though, it is more desirable to think of the wave speed $\sigma$ and elliptic modulus $m$ as ``independent" parameters, that is, the cnoidal solutions are found for fixed elliptic modulus $m \in [0,1]$ and a certain range of wave speed $\sigma>0$.  Indeed, for some cases, it is necessary to assume this condition to have solutions.   For the Schr\"odinger KdV-KdV system \eqref{system1},  these non-trivial periodic traveling-wave solutions are established for each wave speed $\sigma >0$ with $2c>a_0>0$, while for the Schr\"odinger BBM-BBM system \eqref{system2}, $\sigma>0$ with $2c>a_1>0$.  For the Schr\"odinger KdV-BBM \eqref{system3}, the range of wave speed is $\sigma>\frac{a_0}{2c}>0$, while for the Schr\"odinger BBM-KdV \eqref{system4}, $0<\sigma<\frac{2c}{a_1}$.  Moreover, for all four systems, the coefficients $d_2$ and $h_2$ are constant multiples of each other with the ratios being controlled by the coefficients of the third derivatives in the KdV-KdV and BBM-BBM cases as well as the wave speed in the KdV-BBM and BBM-KdV cases.  Precisely, their ratio is an expression of only $a_0$ and $c$ in the KdV-KdV case; $a_1$ and $c$ in the BBM-BBM case; $a_0, c,$ and $\sigma$ in the KdV-BBM case; and $a_1, c,$ and $\sigma$ in the BBM-KdV case.

\section{Exact Jacobi elliptic solutions}
For conciseness, let
\begin{equation}\label{R}
  R = \pm \, \sqrt{m^4-m^2+1}.
\end{equation}  
Then $R\in \mathbb R$ as $m\in[0,1]$.
\subsection{Schr\"odinger KdV-KdV}
Setting all $k_{j,q}=0$ gives us the following set of parameters whenever $2c>a_0 > 0$:
\begin{align*}
  &\begin{cases} 
  B = {\frac {{\it a_0}\,\mu_1-b}{2\,{\it a_0}}} \,,
  \\[5 pt]
  d_1 = h_1 = 0 \,,
  \\[5 pt]
  {\it d_0} = {\frac { \left( {m}^{4}-2\,{m}^{2}R-{m}^{2}+R+1 \right) \sqrt {2\,c-
  {\it a_0}} \left( 3\,{{\it a_0}}^{2}{\mu_1}^{2}-2\,{\it a_0}\,b\mu_1-4\,{
  \it a_0}\,\mu_0-{b}^{2}+4\,{\it a_0} \right) }{8\,\sqrt {{\it a_0}}{R}^{2}
  \left( {\it a_0}-c \right) }} \,,
  \\[5 pt]
  {\it d_2}= {\frac {3\,\sqrt {2\,c-{\it a_0}} \left( 3\,{{\it a_0}}^{2}{
  \mu_1}^{2}-2\,{\it a_0}\,b\mu_1-4\,{\it a_0}\,\mu_0-{b}^{2}+4\,{\it a_0}
  \right) {m}^{2}}{8\,\sqrt {{\it a_0}}R \left( {\it a_0}-c \right) }} \,,
  \\[5 pt]
  {\it h_0}=\frac{-1}{8\,{\it a_0}R\, \left( {\it a_0}-c \right) }\,\Big(6\,{{\it a_0}}^{3}{m}^{2}{\mu_1}^{2}-3\,{{\it a_0}}^{3}{
  \mu_1}^{2}R+6\,{{\it a_0}}^{2}c{\mu_1}^{2}R-3\,{{\it a_0}}^{3}{\mu_1}^{2}-4
  \,{{\it a_0}}^{2}b{m}^{2}\mu_1
  \\
  \quad\quad +2\,{{\it a_0}}^{2}b\mu_1\,R-4\,{\it a_0}\,bc\mu_1\,
  R+2\,{{\it a_0}}^{2}b\mu_1-8\,{{\it a_0}}^{2}{m}^{2}\mu_0+4\,{{\it a_0}}^{2}\mu_0
  \,R-8\,{{\it a_0}}^{2}R\sigma  -2\,{\it a_0}\,{b}^{2}{m}^{2}
  \\
  \quad\quad +{\it a_0}\,{b}^{2}R-
  8\,{\it a_0}\,c\mu_0\,R+8\,{\it a_0}\,cR\sigma-2\,{b}^{2}cR+8\,{{\it a_0}}
  ^{2}{m}^{2}+4\,{{\it a_0}}^{2}\mu_0+4\,{{\it a_0}}^{2}R+{\it a_0}\,{b}^{2}-4\,{{
  \it a_0}}^{2}\Big) \,,
  \end{cases} \displaybreak[3]\\
  &\begin{cases}
  {\it h_2}={\frac { 3\,\left( 3\,{{\it a_0}}^{2}{\mu_1}^{2}-2\,{\it a_0}
  \,b\mu_1-4\,{\it a_0}\,\mu_0-{b}^{2}+4\,{\it a_0} \right) {m}^{2}}{8\,R \left( {
  \it a_0}-c \right) }} \,,
  \\[5 pt]
  \lambda= \sqrt{{\frac {3\,{{\it a_0}}^{2}{\mu_1}^{2}-2\,{\it a_0}\,b\mu_1-4
  \,{\it a_0}\,\mu_0-{b}^{2}+4\,{\it a_0}}{16\,{\it a_0}R\, \left( {\it a_0}-c
  \right) }}} \,,
  \\[5 pt]
  \omega=- \left( {\it a_0}\,{\mu_1}^{2}-\mu_1\,b-\mu_0+\sigma \right) \mu_1 \,,
  \\[5 pt]
  \sigma > 0 \,,
  \\[5 pt]
  m \in [\,0,1\,] \,.
  \end{cases}
\end{align*}
Thus, explicit periodic traveling-wave solutions to the Schr\"odinger KdV-KdV system $\big(u(x,t), v(x,t)\big) = \big(e^{i\omega t} e^{iB(x-\sigma t)}f(x-\sigma t), g(x-\sigma t)\big)$ given in term of the Jacobi cnoidal function 
\begin{equation*}
f(\x)= d_0+d_2 \cn^2(\lambda \xi,m) \quad \text{and} \quad g(\x) = h_0 + h_2 \cn^2(\lambda \xi,m)
\end{equation*}
are established. Notice that $\frac{h_2}{d_2}= 
\sqrt{\frac{a_0}{2c-a_0}}$ and that as $m$ approaches 
$1$, $R$ limits to $\pm 1$.  When $m=R=1$, the above 
coefficients simplify to $d_0=0$ and

$$
\begin{cases}
{\it \tilde h_0}=-\frac{1}{4\,{\it a_0}\, \left( {\it a_0}-c \right) }\,\Big( 3a_0^2c\mu_1^2 -2a_0bc \mu_1 +4a_0^2 - 4a_0c\mu_0 -b^2c -4a_0^2\sigma + 4a_0c \sigma\Big) \,,
\\[5 pt]
{\it \tilde d_2}= {\frac {3\,\sqrt {2\,c-{\it a_0}} \left( 3\,{{\it a_0}}^{2}{
\mu_1}^{2}-2\,{\it a_0}\,b\mu_1-4\,{\it a_0}\,\mu_0-{b}^{2}+4\,{\it a_0}
 \right) }{8\,\sqrt {{\it a_0}} \left( {\it a_0}-c \right) }} \,,
\\[5 pt]
{\it \tilde h_2}={\frac { 3\,\left( 3\,{{\it a_0}}^{2}{\mu_1}^{2}-2\,{\it a_0}
\,b\mu_1-4\,{\it a_0}\,\mu_0-{b}^{2}+4\,{\it a_0} \right) }{8\, \left( {
\it a_0}-c \right) }} \,, \\[5 pt]
\tilde \lambda= \sqrt{{\frac {3\,{{\it a_0}}^{2}{\mu_1}^{2}-2\,{\it a_0}\,b\mu_1-4
\,{\it a_0}\,\mu_0-{b}^{2}+4\,{\it a_0}}{16\,{\it a_0}\, \left( {\it a_0}-c
 \right) }}} \,,
\\[5 pt]
\omega=- \left( {\it a_0}\,{\mu_1}^{2}-\mu_1\,b-\mu_0+\sigma \right) \mu_1 \,,
\end{cases}
$$
from which one obtains the following solitary-wave solution to the system \eqref{system1} 
\begin{equation*}
u(x,t) = e^{i\omega t} e^{iB(x-\sigma t)} \tilde f(x-\sigma t) \ \ \ \ \ \mbox{and} \ \ \ \ \ v(x,t) =\tilde h_0  \pm\sqrt{ \dd\frac{a_0}{2c-a_0}} \tilde f(x-\sigma t)
\end{equation*}
where $\tilde f(\xi) = \tilde d_2 \operatorname{sech}^2 (\tilde \lambda \xi) $.  Furthermore, when $\sigma = \frac{4a_0^2 + 3a_0^2 c \mu_1^2 -2 a_0bc \mu_1 - 4a_0c\mu_0 -b^2c}{4a_0(a_0-c)}$, one has $\tilde h_0=0$ and the synchronized solitary-wave solution established in \cite{NLB} is recovered.  

When $m=-R=1$, the above coefficients simplify to
$$
\begin{cases}
{\it \bar d_0} = {\frac {\sqrt {2\,c-
{\it a_0}} \left( 3\,{{\it a_0}}^{2}{\mu_1}^{2}-2\,{\it a_0}\,b\mu_1-4\,{
\it a_0}\,\mu_0-{b}^{2}+4\,{\it a_0} \right) }{4\sqrt {{\it a_0}}
 \left( {\it a_0}-c \right) }} \,,
\\[5 pt]
{\it \bar h_0}= \frac{3a_0^3\mu_1^2 -2a_0^2b \mu_1 -4a_0^2 \mu_0 -a_0b^2 - 3 a_0^2c \mu_1^2 + 2a_0bc\mu_1 + 4a_0c \mu_0 +b^2 c + 4a_0 \sigma (a_0 -c)}{ 4\,{\it a_0}\, \left( {\it a_0}-c \right) }\,,
\\[5 pt]
{\it \bar d_2}= {\frac {3\,\sqrt {2\,c-{\it a_0}} \left( 3\,{{\it a_0}}^{2}{
\mu_1}^{2}-2\,{\it a_0}\,b\mu_1-4\,{\it a_0}\,\mu_0-{b}^{2}+4\,{\it a_0}
 \right) }{8\,\sqrt {{\it a_0}} \left( {\it a_0}-c \right) }} \,,
\\[5 pt]
{\it \bar h_2}=-{\frac { 3\,\left( 3\,{{\it a_0}}^{2}{\mu_1}^{2}-2\,{\it a_0}
\,b\mu_1-4\,{\it a_0}\,\mu_0-{b}^{2}+4\,{\it a_0} \right) }{8\, \left( {
\it a_0}-c \right) }} \,,
\\[5 pt]
\bar \lambda= \sqrt{{-\frac {3\,{{\it a_0}}^{2}{\mu_1}^{2}-2\,{\it a_0}\,b\mu_1-4
\,{\it a_0}\,\mu_0-{b}^{2}+4\,{\it a_0}}{16\,{\it a_0}\, \left( {\it a_0}-c
 \right) }}} \,,
\\[5 pt]
\omega=- \left( {\it a_0}\,{\mu_1}^{2}-\mu_1\,b-\mu_0+\sigma \right) \mu_1 \,,
\end{cases}
$$
and one arrives at the solitary-wave solution 
\begin{equation*}
u(x,t) = e^{i\omega t} e^{iB(x-\sigma t)}[ \bar d_0 + \bar f(x-\sigma t)] \ \ \ \ \ \mbox{and} \ \ \ \ \ v(x,t) =\bar h_0  \pm \sqrt{\dd\frac{a_0}{2c-a_0}}\bar f(x-\sigma t)
\end{equation*} where $\bar f(\xi) = \bar d_2 \operatorname{sech}^2 (\bar \lambda \xi) $.

Aside from the above non-trivial solutions, system \eqref{system1} also possessess the following trivial and semi-trivial solutions
\begin{enumerate}
\item
\begin{equation*}
u(x,t) =0  \ \ \ \ \ \mbox{and} \ \ \ \ \ v(x,t) =h_0;
\end{equation*} for any $h_0 \in \mathbb R$.
\item
\begin{equation*}
u(x,t) = e^{i\omega t} e^{iB(x-\sigma t)} d_0  \ \ \ \ \ \mbox{and} \ \ \ \ \ v(x,t) =h_0;
\end{equation*} where $\s = \frac{\omega - a_0B^3 - b B^2 + B h_0 + B \mu_0 + h_0 \mu_1}{B}$, for any $B, d_0, h_0, \omega \in \mathbb R$. 
\item
\begin{equation*}
u(x,t) = 0  \ \ \ \ \ \mbox{and} \ \ \ \ \ v(x,t) = -\frac{2}{3}h_2+\frac{1}{3}\,{\frac {h_2}{{m}^{2}}}+\sigma-1 +h_2\cn^2\big(\lambda (x-\sigma t),m\big) ;
\end{equation*} where $\lambda= \sqrt{\frac{h_2}{12c{m}^{2}}} $, for any $h_2, \sigma>0$ and $m \in [\,0,1\,]$.
\end{enumerate}

\subsection{Schr\"odinger BBM-BBM}
Setting all $k_{j,q}=0$ gives us the following whenever $2c>a_1>0$ and $R$ is as defined in \eqref{R}:

\begin{align*}
&\begin{cases}
B = {\frac {{\it a_1}\,\mu_0\,\mu_1-b}{2a_1\sigma\, \left( {\it a_1
}\,{\mu_1}^{2}+1 \right) }} \,,
\\[5 pt]
d_1 = h_1 = 0 \,,
\\[5 pt]
{\it d_0}= \frac{\sqrt {{\it a_1}\, \left( 2\,c - {\it a_1} \right) 
} \left( {m}^{4}+2\,{m}^{2}R-{m}^{2}-R+1 \right)}{8\,{\it a_1}{R}^{2}
\sigma\, \left( {\it a_1}\,{\mu_1}^{2}+1 \right) ^{2} \left( {
\it a_1}-c \right) } \Big( 4\,{{\it a_1}}^{3}{\mu_1}^{4}\sigma-4\,{{\it a_1}}^{2}b{\mu_1}^{3
}\sigma-{{\it a_1}}^{2}{\mu_0}^{2}{\mu_1}^{2}-4\,{{\it a_1}}^{2}\mu_0\,{\mu
1}^{2}\sigma+8\,{{\it a_1}}^{2}{\mu_1}^{2}\sigma
\\
\quad\quad +2\,{\it a_1}\,b\mu_0\,\mu
1-4\,{\it a_1}\,b\mu_1\,\sigma-4\,{\it a_1}\,\mu_0\,\sigma+4\,{\it a_1}\,
\sigma-{b}^{2} \Big)  \,,
\\[5 pt]
{\it d_2}=\frac{-3\,{m}^{2} \sqrt {{\it a_1}\, \left( 2\,c - {\it a_1}
 \right) }}{8\,a_1R\sigma\, \left( {\it a_1}
\,{\mu_1}^{2}+1 \right) ^{2} \left( {\it a_1}-c \right)} \Big( 4\,{{\it a_1}}^{3}{\mu_1}^{4}\sigma-4\,{{\it a_1}}^{2}b
{\mu_1}^{3}\sigma-{{\it a_1}}^{2}{\mu_0}^{2}{\mu_1}^{2}-4\,{{\it a_1}}^{2}
\mu_0\,{\mu_1}^{2}\sigma+8\,{{\it a_1}}^{2}{\mu_1}^{2}\sigma
\\
\quad\quad +2\,{\it a_1}\,
b\mu_0\,\mu_1-4\,{\it a_1}\,b\mu_1\,\sigma-4\,{\it a_1}\,\mu_0\,\sigma+4\,{
\it a_1}\,\sigma-{b}^{2} \Big) \,,
\\[5 pt]
{\it h_0}=\frac{1}{8\,a_1R\sigma\,
 \left( {\it a_1}\,{\mu_1}^{2}+1 \right) ^{2} \left( {\it a_1}-c \right)}\Big(8\,{{\it a_1}}^{4}{\mu_1}^{4}r{\sigma}^{2}-8\,{{
\it a_1}}^{3}c{\mu_1}^{4}R{\sigma}^{2}+8\,{{\it a_1}}^{4}{m}^{2}{\mu_1}^{4}
\sigma-4\,{{\it a_1}}^{4}{\mu_1}^{4}R\sigma-4\,{{\it a_1}}^{4}{\mu_1}^{4}
\sigma
\\
\quad\quad -8\,{{\it a_1}}^{3}b{m}^{2}{\mu_1}^{3}\sigma-4\,{{\it a_1}}^{3}b{\mu_1}^{3
}R\sigma+8\,{{\it a_1}}^{2}bc{\mu_1}^{3}R\sigma+4\,{{\it a_1}}^{3}b{\mu_1}
^{3}\sigma-2\,{{\it a_1}}^{3}{m}^{2}{\mu_0}^{2}{\mu_1}^{2}-8\,{{\it a_1}}^{3}{m}^{2}
\mu_0\,{\mu_1}^{2}\sigma
\\
\quad\quad -{{\it a_1}}^{3}{\mu_0}^{2}{\mu_1}^{2}R-4\,{{\it a_1
}}^{3}\mu_0\,{\mu_1}^{2}R\sigma+16\,{{\it a_1}}^{3}{\mu_1}^{2}R{\sigma}^{2
}+2\,{{\it a_1}}^{2}c{\mu_0}^{2}{\mu_1}^{2}R+8\,{{\it a_1}}^{2}c\mu_0\,{\mu
1}^{2}R\sigma-16\,{{\it a_1}}^{2}c{\mu_1}^{2}R{\sigma}^{2}
\\
\quad\quad +16\,{{\it a_1}
}^{3}{m}^{2}{\mu_1}^{2}\sigma+{{\it a_1}}^{3}{\mu_0}^{2}{\mu_1}^{2}+4\,{{\it a_1}
}^{3}\mu_0\,{\mu_1}^{2}\sigma-8\,{{\it a_1}}^{3}{\mu_1}^{2}R\sigma-8\,{{
\it a_1}}^{3}{\mu_1}^{2}\sigma+4\,{{\it a_1}}^{2}b{m}^{2}\mu_0\,\mu_1
\\
\quad\quad -8\,{{\it a_1
}}^{2}b{m}^{2}\mu_1\,\sigma+2\,{{\it a_1}}^{2}b\mu_0\,\mu_1\,R-4\,{{\it a_1}}^{2}
b\mu_1\,R\sigma-4\,{\it a_1}\,bc\mu_0\,\mu_1\,R+8\,{\it a_1}\,bc\mu_1\,R
\sigma-2\,{{\it a_1}}^{2}\,b\,\mu_0\,\mu_1
\\
\quad\quad+4\,{{\it a_1}}^{2}b\mu_1\,\sigma-8
\,{{\it a_1}}^{2}{m}^{2}\mu_0\,\sigma-4\,{{\it a_1}}^{2}\mu_0\,R\sigma+8\,{{\it 
a_1}}^{2}R{\sigma}^{2}+8\,{\it a_1}\,c\mu_0\,R\sigma-8\,{\it a_1}\,cR{
\sigma}^{2}+8\,{{\it a_1}}^{2}{m}^{2}\sigma
\\
\quad\quad +4\,{{\it a_1}}^{2}\mu_0\,\sigma-4\,
{{\it a_1}}^{2}R\sigma-2\,{\it a_1}\,{b}^{2}{m}^{2}-{\it a_1}\,{b}^{2}R+2\,{b}^
{2}cR-4\,{{\it a_1}}^{2}\sigma+{\it a_1}\,{b}^{2}\Big) \,,
\\[5 pt]
{\it h_2}={\frac { -3\,\left( 4\,{{\it a_1}}^{3}{\mu_1}^{4}\sigma-4\,{{
\it a_1}}^{2}b{\mu_1}^{3}\sigma-{{\it a_1}}^{2}{\mu_0}^{2}{\mu_1}^{2}-4\,{{
\it a_1}}^{2}\mu_0\,{\mu_1}^{2}\sigma+8\,{{\it a_1}}^{2}{\mu_1}^{2}\sigma+2
\,{\it a_1}\,b\mu_0\,\mu_1-4\,{\it a_1}\,b\mu_1\,\sigma-4\,{\it a_1}\,\mu_0\,
\sigma+4\,{\it a_1}\,\sigma-{b}^{2} \right) {m}^{2}}{8\,R\sigma\, \left( {\it a_1
}\,{\mu_1}^{2}+1 \right) ^{2} \left( {\it a_1}-c \right) }} \,,
\end{cases} \displaybreak[3]\\
&\begin{cases}
\lambda=\sqrt{\frac {4\,{{\it a_1}}^{3}{\mu_1}^{4}\sigma-4\,{{\it a_1}}
^{2}b{\mu_1}^{3}\sigma-{{\it a_1}}^{2}{\mu_0}^{2}{\mu_1}^{2}-4\,{{\it a_1}}
^{2}\mu_0\,{\mu_1}^{2}\sigma+8\,{{\it a_1}}^{2}{\mu_1}^{2}\sigma+2\,{\it 
a_1}\,b\mu_0\,\mu_1-4\,{\it a_1}\,b\mu_1\,\sigma-4\,{\it a_1}\,\mu_0\,\sigma+
4\,{\it a_1}\,\sigma-{b}^{2}}{ -16\,a_1R\sigma^2\left( {\it a_1}-c \right)
 \left( {\it a_1}\,{\mu_1}^{2}+1 \right) ^{2}}} \,,
\\[5 pt]
\omega=-{\frac { \left( {\it a_1}\,{\mu_1}^{2}\sigma-b\mu_1-\mu_0+\sigma
 \right) \mu_1}{{\it a_1}\,{\mu_1}^{2}+1}} \,,
\\[5 pt]
\sigma > 0  \,,
\\[5 pt]
m \in [\,0,1\,] \,.
\end{cases}
\end{align*}
Thus, explicit periodic traveling-wave solutions to the Schr\"odinger BBM-BBM system $\big(u(x,t), v(x,t)\big) = \big(e^{i\omega t} e^{iB(x-\sigma t)}f(x-\sigma t), g(x-\sigma t)\big)$  given in term of the Jacobi cnoidal function
\begin{equation*}
f(\x)= d_0+d_2 \cn^2(\lambda \xi,m) \quad \text{and} \quad g(\x) = h_0 + h_2 \cn^2(\lambda \xi,m)
\end{equation*} are established.  Notice that $\frac{h_2}{d_2}= \sqrt{\frac{a_1}{2c-a_1}}$.   When $m=R=1$, the above coefficients simplify to
$$
\begin{cases}
{\it \tilde d_0}= \frac{\sqrt {{\it a_1}\, \left( 2\,c - {\it a_1} \right) 
} }{4 a_1\sigma\, \left( {\it a_1}\,{\mu_1}^{2}+1 \right) ^{2} \left( {
\it a_1}-c \right) } \Big( 4\,{{\it a_1}}^{3}{\mu_1}^{4}\sigma-4\,{{\it a_1}}^{2}b{\mu_1}^{3
}\sigma-{{\it a_1}}^{2}{\mu_0}^{2}{\mu_1}^{2}-4\,{{\it a_1}}^{2}\mu_0\,{\mu_1}^{2}\sigma+8\,{{\it a_1}}^{2}{\mu_1}^{2}\sigma
\\
\quad\quad +2\,{\it a_1}\,b\mu_0\,\mu_1-4\,{\it a_1}\,b\mu_1\,\sigma-4\,{\it a_1}\,\mu_0\,\sigma+4\,{\it a_1}\,
\sigma-{b}^{2} \Big)  \,,
\\[5 pt]
{\it \tilde d_2}=\frac{-3 \sqrt {{\it a_1}\, \left( 2\,c - {\it a_1}
 \right) }}{8a_1\sigma\, \left( {\it a_1}
\,{\mu_1}^{2}+1 \right) ^{2} \left( {\it a_1}-c \right)} \Big( 4\,{{\it a_1}}^{3}{\mu_1}^{4}\sigma-4\,{{\it a_1}}^{2}b
{\mu_1}^{3}\sigma-{{\it a_1}}^{2}{\mu_0}^{2}{\mu_1}^{2}-4\,{{\it a_1}}^{2}
\mu_0\,{\mu_1}^{2}\sigma+8\,{{\it a_1}}^{2}{\mu_1}^{2}\sigma
\\
\quad\quad +2\,{\it a_1}\,
b\mu_0\,\mu_1-4\,{\it a_1}\,b\mu_1\,\sigma-4\,{\it a_1}\,\mu_0\,\sigma+4\,{
\it a_1}\,\sigma-{b}^{2} \Big) \,,
\\[5 pt]
{\it \tilde h_0}=\frac{1}{4a_1\sigma\,
 \left( {\it a_1}\,{\mu_1}^{2}+1 \right) ^{2} \left( {\it a_1}-c \right)}\Big(4\,{{\it a_1}}^{4}{\mu_1}^{4}{\sigma}^{2}-4\,{{
\it a_1}}^{3}c{\mu_1}^{4}{\sigma}^{2} -4\,{{\it a_1}}^{3}b{\mu_1}^{3}\sigma+4\,{{\it a_1}}^{2}bc{\mu_1}^{3}\sigma-\,{{\it a_1}}^{3}{\mu_0}^{2}{\mu_1}^{2}\\
\quad\quad-4\,{{\it a_1}}^{3}
\mu_0\,{\mu_1}^{2}\sigma +8\,{{\it a_1}}^{3}{\mu_1}^{2}{\sigma}^{2
}+\,{{\it a_1}}^{2}c{\mu_0}^{2}{\mu_1}^{2}+4\,{{\it a_1}}^{2}c\mu_0\,{\mu_1}^{2}\sigma-8\,{{\it a_1}}^{2}c{\mu_1}^{2}{\sigma}^{2}
+2\,{{\it a_1}}^{2}b\mu_0\,\mu_1 \\
\quad\quad-4\,{{\it a_1
}}^{2}b\mu_1\,\sigma-2\,{\it a_1}\,bc\mu_0\,\mu_1+4\,{\it a_1}\,bc\mu_1\sigma
-4\,{{\it a_1}}^{2}\mu_0\,\sigma+4\,{{\it 
a_1}}^{2}{\sigma}^{2}+4\,{\it a_1}\,c\mu_0\sigma-4\,{\it a_1}\,c{
\sigma}^{2}\\
\quad\quad
-\,{\it a_1}\,{b}^{2}+\,{b}^
{2}c\Big) \,,
\\[5 pt]
{\it \tilde h_2}={\frac { -3\,\left( 4\,{{\it a_1}}^{3}{\mu_1}^{4}\sigma-4\,{{
\it a_1}}^{2}b{\mu_1}^{3}\sigma-{{\it a_1}}^{2}{\mu_0}^{2}{\mu_1}^{2}-4\,{{
\it a_1}}^{2}\mu_0\,{\mu_1}^{2}\sigma+8\,{{\it a_1}}^{2}{\mu_1}^{2}\sigma+2
\,{\it a_1}\,b\mu_0\,\mu_1-4\,{\it a_1}\,b\mu_1\,\sigma-4\,{\it a_1}\,\mu_0\,
\sigma+4\,{\it a_1}\,\sigma-{b}^{2} \right) }{8\sigma\, \left( {\it a_1
}\,{\mu_1}^{2}+1 \right) ^{2} \left( {\it a_1}-c \right) }} \,,
\\[5 pt]
\tilde \lambda=\sqrt{\frac {4\,{{\it a_1}}^{3}{\mu_1}^{4}\sigma-4\,{{\it a_1}}
^{2}b{\mu_1}^{3}\sigma-{{\it a_1}}^{2}{\mu_0}^{2}{\mu_1}^{2}-4\,{{\it a_1}}
^{2}\mu_0\,{\mu_1}^{2}\sigma+8\,{{\it a_1}}^{2}{\mu_1}^{2}\sigma+2\,{\it 
a_1}\,b\mu_0\,\mu_1-4\,{\it a_1}\,b\mu_1\,\sigma-4\,{\it a_1}\,\mu_0\,\sigma+
4\,{\it a_1}\,\sigma-{b}^{2}}{ -16a_1\sigma^2\,\left( {\it a_1}-c \right) 
 \left( {\it a_1}\,{\mu_1}^{2}+1 \right) ^{2}}} \,,
\\[5 pt]
\omega=-{\frac { \left( {\it a_1}\,{\mu_1}^{2}\sigma-b\mu_1-\mu_0+\sigma
 \right) \mu_1}{{\it a_1}\,{\mu_1}^{2}+1}} \,,
\\[5 pt]
\sigma > 0  \,,
\end{cases}
$$
from which one obtains the following solitary-wave solution to the system \eqref{system2} 
\begin{equation*}
u(x,t) = e^{i\omega t} e^{iB(x-\sigma t)}[\tilde d_0 + \tilde f(x-\sigma t)] \ \ \ \ \ \mbox{and} \ \ \ \ \ v(x,t) =\tilde h_0  \pm \sqrt{\dd\frac{a_1}{2c-a_1}} \tilde f(x-\sigma t)
\end{equation*}
where $\tilde f(\xi) = \tilde d_2 \operatorname{sech}^2 (\tilde \lambda \xi) $. 

When $m=-R=1$, the above coefficients simplify to $d_0=0$ and

\begin{align*}
&\begin{cases}
{\it \bar d_2}=\frac{3 \sqrt {{\it a_1}\, \left( 2\,c - {\it a_1}
 \right) }}{8a_1\sigma\, \left( {\it a_1}
\,{\mu_1}^{2}+1 \right) ^{2} \left( {\it a_1}-c \right)} \Big( 4\,{{\it a_1}}^{3}{\mu_1}^{4}\sigma-4\,{{\it a_1}}^{2}b
{\mu_1}^{3}\sigma-{{\it a_1}}^{2}{\mu_0}^{2}{\mu_1}^{2}-4\,{{\it a_1}}^{2}
\mu_0\,{\mu_1}^{2}\sigma+8\,{{\it a_1}}^{2}{\mu_1}^{2}\sigma
\\
\quad\quad +2\,{\it a_1}\,
b\mu_0\,\mu_1-4\,{\it a_1}\,b\mu_1\,\sigma-4\,{\it a_1}\,\mu_0\,\sigma+4\,{
\it a_1}\,\sigma-{b}^{2} \Big) \,,
\\[5 pt]
{\it \bar h_0}=\frac{1}{-8a_1\sigma\,
 \left( {\it a_1}\,{\mu_1}^{2}+1 \right) ^{2} \left( {\it a_1}-c \right)}\Big(-8\,{{\it a_1}}^{4}{\mu_1}^{4}{\sigma}^{2}+8\,{{
\it a_1}}^{3}c{\mu_1}^{4}{\sigma}^{2}+8\,{{\it a_1}}^{4}{\mu_1}^{4}
\sigma -8\,{{\it a_1}}^{2}bc{\mu_1}^{3}\sigma
-16\,{{\it a_1}}^{3}{\mu_1}^{2}{\sigma}^{2
}\\
\quad\quad -2\,{{\it a_1}}^{2}c{\mu_0}^{2}{\mu_1}^{2}-8\,{{\it a_1}}^{2}c\mu_0\,{\mu_1}^{2}\sigma+16\,{{\it a_1}}^{2}c{\mu_1}^{2}{\sigma}^{2}
 +16\,{{\it a_1}
}^{3}{\mu_1}^{2}\sigma+4\,{\it a_1}\,bc\mu_0\,\mu_1-8\,{\it a_1}\,bc\mu_1\sigma
\end{cases} \displaybreak[3]\\
&\begin{cases}
\quad\quad-8\,{{\it 
a_1}}^{2}{\sigma}^{2}-8\,{\it a_1}\,c\mu_0\sigma+8\,{\it a_1}\,c{
\sigma}^{2}+8\,{{\it a_1}}^{2}\sigma -2\,{b}^
{2}c\Big) \,,
\\[5 pt]
{\it \bar h_2}={\frac {3\,\left( 4\,{{\it a_1}}^{3}{\mu_1}^{4}\sigma-4\,{{
\it a_1}}^{2}b{\mu_1}^{3}\sigma-{{\it a_1}}^{2}{\mu_0}^{2}{\mu_1}^{2}-4\,{{
\it a_1}}^{2}\mu_0\,{\mu_1}^{2}\sigma+8\,{{\it a_1}}^{2}{\mu_1}^{2}\sigma+2
\,{\it a_1}\,b\mu_0\,\mu_1-4\,{\it a_1}\,b\mu_1\,\sigma-4\,{\it a_1}\,\mu_0\,
\sigma+4\,{\it a_1}\,\sigma-{b}^{2} \right)}{8\sigma\, \left( {\it a_1
}\,{\mu_1}^{2}+1 \right) ^{2} \left( {\it a_1}-c \right) }} \,,
\\[5 pt]
\bar \lambda=\sqrt{\frac {4\,{{\it a_1}}^{3}{\mu_1}^{4}\sigma-4\,{{\it a_1}}
^{2}b{\mu_1}^{3}\sigma-{{\it a_1}}^{2}{\mu_0}^{2}{\mu_1}^{2}-4\,{{\it a_1}}
^{2}\mu_0\,{\mu_1}^{2}\sigma+8\,{{\it a_1}}^{2}{\mu_1}^{2}\sigma+2\,{\it 
a_1}\,b\mu_0\,\mu_1-4\,{\it a_1}\,b\mu_1\,\sigma-4\,{\it a_1}\,\mu_0\,\sigma+
4\,{\it a_1}\,\sigma-{b}^{2}}{16a_1\sigma^2\,\left( {\it a_1}-c \right)
 \left( {\it a_1}\,{\mu_1}^{2}+1 \right) ^{2}}} \,,
\\[5 pt]
\omega=-{\frac { \left( {\it a_1}\,{\mu_1}^{2}\sigma-b\mu_1-\mu_0+\sigma
 \right) \mu_1}{{\it a_1}\,{\mu_1}^{2}+1}} \,,
\end{cases}
\end{align*}
and one arrives at the solitary-wave solution 
\begin{equation*}
u(x,t) = e^{i\omega t} e^{iB(x-\sigma t)} \bar f(x-\sigma t) \ \ \ \ \ \mbox{and} \ \ \ \ \ v(x,t) =\bar h_0  \pm \sqrt{\dd\frac{a_1}{2c-a_1}}\bar f(x-\sigma t)
\end{equation*} where $\bar f(\xi) = \bar d_2 \operatorname{sech}^2 (\bar \lambda \xi) $.  Furthermore, when $ B = {\frac {{\it a_1}\,\mu_0\,\mu_1-b}{2\,\sigma\,{\it a_1}\, \left( {\it a_1
}\,{\mu_1}^{2}+1 \right) }} $ satisfies the following equation
\begin{equation*}
    (a_1^2 c \mu_0 \mu_1 - a_1 b c) B^2 + (2 a_1 b c \mu_1 + 2 a_1 c \mu_0 - 2 a_1^3 \mu_1^2 - 2 a_1^2) B + (a_1^2 \mu_0 \mu_1 + b c - a_1 b - a_1 c \mu_0 \mu_1) = 0,
\end{equation*}
one has $\bar h_0 =0$ and the synchronized solitary-wave solution established in \cite{NLB} is recovered.  

Aside from the above non-trivial solutions, system \eqref{system2} also possessess the following trivial and semi-trivial solutions 
\begin{enumerate}
\item
\begin{equation*}
u(x,t) =0  \ \ \ \ \ \mbox{and} \ \ \ \ \ v(x,t) =h_0;
\end{equation*} for any $h_0 \in \mathbb R$.
\item
\begin{equation*}
u(x,t) = e^{i\omega t} e^{iB(x-\sigma t)} d_0  \ \ \ \ \ \mbox{and} \ \ \ \ \ v(x,t) =h_0;
\end{equation*} where $\s = \frac{a_1 B^2 \omega - b B^2 + B h_0 + B \mu_0 + h_0 \mu_1 + \omega}{B (a_1 B^2 + 1)} $, for any $B, d_0, h_0,\omega \in \mathbb R$.
\item
\begin{equation*}
u(x,t) = 0  \ \ \ \ \ \mbox{and} \ \ \ \ \ v(x,t) = -\frac{2}{3}h_2+\frac{1}{3}\,{\frac {h_2}{{m}^{2}}}+\sigma-1 + h_2\cn^2\big(\lambda (x-\sigma t),m\big) ;
\end{equation*} where $ \l = \sqrt{\frac{h_2}{12c{m}^{2}\sigma}}$, for any $h_2, \sigma>0$ and $m \in [\,0,1\,]$. 
\end{enumerate}

\subsection{Schr\"odinger KdV-BBM}
Setting all $k_{j,q}=0$ gives us the following set of parameters with $R$ as defined in \eqref{R}:

\begin{align*}
&\begin{cases}
B = {\frac {{\it a_0}\,\mu_1-b}{2\,{\it a_0}}} \,,
\\[5 pt]
d_1 = h_1 = 0 \,,
\\[5 pt]
{\it d_0}= {\frac { \left( {m}^{4}-2\,{m}^{2}R-{m}^{2}+R+1 \right) \sqrt {2\,c
\sigma-{\it a_0}} \left( 3\,{{\it a_0}}^{2}{\mu_1}^{2}-2\,{\it a_0}\,b\mu_1
-4\,{\it a_0}\,\mu_0-{b}^{2}+4\,{\it a_0} \right) }{8\sqrt {{\it a_0
}}\,{R}^{2} \left( {\it a_0} -c\sigma \right) }} \,,
\\[5 pt]
{\it d_2}= {\frac {3\,\sqrt {2\,c\sigma-{\it a_0}} \left( 3\,{{\it a_0}}
^{2}{\mu_1}^{2}-2\,{\it a_0}\,b\mu_1-4\,{\it a_0}\,\mu_0-{b}^{2}+4\,{\it a_0
} \right) {m}^{2}}{8\,\sqrt {{\it a_0}}R \left( {\it a_0} -c\sigma \right) }} \,,
\\[5 pt]
{\it h_0}= \frac{1}{8{\it a_0}R
\, \left(c\sigma-{\it a_0} \right) }\,\Big( 6\,{{\it a_0}}^{2}c{\mu_1}^{2}R\sigma+6\,{{\it 
a_0}}^{3}{m}^{2}{\mu_1}^{2}-3\,{{\it a_0}}^{3}{\mu_1}^{2}R-4\,{\it a_0}\,bc\mu_1\,
R\sigma
\\
\quad\quad -3\,{{\it a_0}}^{3}{\mu_1}^{2}-4\,{{\it a_0}}^{2}b{m}^{2}\mu_1+2\,{{
\it a_0}}^{2}b\mu_1\,R-8\,{\it a_0}\,c\mu_0\,R\sigma+8\,{\it a_0}\,cR{
\sigma}^{2}-2\,{b}^{2}cR\sigma
\\
\quad\quad +2\,{{\it a_0}}^{2}b\mu_1-8\,{{\it a_0}}
^{2}{m}^{2}\mu_0+4\,{{\it a_0}}^{2}\mu_0\,R-8\,{{\it a_0}}^{2}R\sigma-2\,{\it 
a_0}\,{b}^{2}{m}^{2}+{\it a_0}\,{b}^{2}R
\end{cases} \displaybreak[3]\\
&\begin{cases}
\quad\quad +8\,{{\it a_0}}^{2}{m}^{2}+4\,{{\it a_0}}^{2}
\mu_0+4\,{{\it a_0}}^{2}R+{\it a_0}\,{b}^{2}-4\,{{\it a_0}}^{2}\Big) \,,
\\[5 pt]
{\it h_2}= {\frac {3\, \left( 3\,{{\it a_0}}^{2}{\mu_1}^{2}-2\,{\it a_0}
\,b\mu_1-4\,{\it a_0}\,\mu_0-{b}^{2}+4\,{\it a_0} \right) {m}^{2}}{8R \left( {\it a_0} 
-c\sigma \right) }} \,,
\\[5 pt]
\lambda= \sqrt{{\frac {3\,{{\it a_0}}^{2}{\mu_1}^{2}-2\,{\it a_0}\,b\mu_1-4
\,{\it a_0}\,\mu_0-{b}^{2}+4\,{\it a_0}}{16{\it a_0}R\, \left( {\it a_0} 
-c\sigma \right) }}}
 \,,
\\[5 pt]
\omega=- \left( {\it a_0}\,{\mu_1}^{2}-\mu_1\,b-\mu_0+\sigma \right) \mu_1 \,,
\\[5 pt]
\sigma > \frac{a_0}{2\,c} \,,
\\[5 pt]
m \in [\,0,1\,] \,.
\end{cases}
\end{align*}
Thus, explicit periodic traveling-wave solutions to the Schr\"odinger KdV-BBM system $\big(u(x,t), v(x,t)\big) = \big(e^{i\omega t} e^{iB(x-\sigma t)}f(x-\sigma t), g(x-\sigma t)\big)$  given in term of the Jacobi cnoidal function
\begin{equation*}
f(\x)= d_0+d_2 \cn^2(\lambda \xi,m) \quad \text{and} \quad g(\x) = h_0 + h_2 \cn^2(\lambda \xi,m)
\end{equation*} 
are established. Notice that $\frac{h_2}{d_2}= \sqrt{\frac{a_0}{2c\sigma-a_0}}$.   When $m=R=1$, the above coefficients simplify to $d_0=0$ and

$$
\begin{cases}
{\it \tilde d_2}= {\frac {3\,\sqrt {2\,c\sigma-{\it a_0}} \left( 3\,{{\it a_0}}
^{2}{\mu_1}^{2}-2\,{\it a_0}\,b\mu_1-4\,{\it a_0}\,\mu_0-{b}^{2}+4\,{\it a_0
} \right) }{8\,\sqrt {{\it a_0}} \left( {\it a_0} -c\sigma \right) }} \,,
\\[5 pt]
{\it \tilde h_0}= \frac{1}{8{\it a_0}
\, \left(c\sigma-{\it a_0} \right) }\,\Big( 6\,{{\it a_0}}^{2}c{\mu_1}^{2}\sigma-4\,{\it a_0}\,bc\mu_1\,
\sigma -8\,{\it a_0}\,c\mu_0\,\sigma+8\,{\it a_0}\,c{
\sigma}^{2}-2\,{b}^{2}c\sigma-8\,{{\it a_0}}^{2}\sigma+8\,{{\it a_0}}^{2}\Big) \,,
\\[5 pt]
{\it \tilde h_2}= {\frac {3\, \left( 3\,{{\it a_0}}^{2}{\mu_1}^{2}-2\,{\it a_0}
\,b\mu_1-4\,{\it a_0}\,\mu_0-{b}^{2}+4\,{\it a_0} \right)}{8 \left( {\it a_0} 
-c\sigma \right) }} \,,
\\[5 pt]
\tilde \lambda= \sqrt{{\frac {3\,{{\it a_0}}^{2}{\mu_1}^{2}-2\,{\it a_0}\,b\mu_1-4
\,{\it a_0}\,\mu_0-{b}^{2}+4\,{\it a_0}}{16{\it a_0}\, \left( {\it a_0} -c\sigma \right) }}}
 \,,
\\[5 pt]
\omega=- \left( {\it a_0}\,{\mu_1}^{2}-\mu_1\,b-\mu_0+\sigma \right) \mu_1 \,,
\\[5 pt]
\sigma > \frac{a_0}{2\,c} \,,
\end{cases}
$$
from which one obtains the following solitary-wave solution to the system \eqref{system3} 
\begin{equation*}
u(x,t) = e^{i\omega t} e^{iB(x-\sigma t)} \tilde f(x-\sigma t) \ \ \ \ \ \mbox{and} \ \ \ \ \ v(x,t) =\tilde h_0  \pm\sqrt{ \dd\frac{a_0}{2c\sigma-a_0}} \tilde f(x-\sigma t)
\end{equation*}
where $\tilde f(\xi) = \tilde d_2 \operatorname{sech}^2 (\tilde \lambda \xi) $.   Furthermore, when $\sigma$ satisfies the condition
$ \frac{\s + 3 a_0 B^2 + 2 b B - \mu_0}{a_0} = \frac{\s - 1}{c \s}, $
one has $\tilde h_0=0$ and the synchronized solitary-wave solution established in \cite{NLB} is recovered.  

When $m=-R=1$, the above coefficients simplify to
\begin{align*}
&\begin{cases}
{\it \bar d_0}= {\frac { \sqrt {2\,c
\sigma-{\it a_0}} \left( 3\,{{\it a_0}}^{2}{\mu_1}^{2}-2\,{\it a_0}\,b\mu_1
-4\,{\it a_0}\,\mu_0-{b}^{2}+4\,{\it a_0} \right) }{4\sqrt {{\it a_0
}} \left( {\it a_0} -c\sigma \right) }} \,,
\\[5 pt]
{\it \bar d_2}= {\frac {3\,\sqrt {2\,c\sigma-{\it a_0}} \left( 3\,{{\it a_0}}
^{2}{\mu_1}^{2}-2\,{\it a_0}\,b\mu_1-4\,{\it a_0}\,\mu_0-{b}^{2}+4\,{\it a_0
} \right) }{8\,\sqrt {{\it a_0}}\left( c\sigma - a_0 \right) }} \,,
\\[5 pt]
{\it \bar h_0}= \frac{1}{4{\it a_0}
\, \left(c\sigma-{\it a_0} \right) }\,\Big( 3\,{{\it a_0}}^{2}c{\mu_1}^{2}\sigma
\, -3{{\it 
a_0}}^{3}{\mu_1}^{2}-2\,{\it a_0}\,bc\mu_1\,
\sigma
+2\,{{
\it a_0}}^{2}b\mu_1-4\,{\it a_0}\,c\mu_0\,\sigma+4\,{\it a_0}\,c{
\sigma}^{2}\\
\quad\quad -\,{b}^{2}c\sigma
+4\,{{\it a_0}}
^{2}\mu_0-4\,{{\it a_0}}^{2}\sigma+{\it a_0}\,{b}^{2}\Big)\,,
\\[5 pt]
{\it \bar h_2}= {\frac {3\, \left( 3\,{{\it a_0}}^{2}{\mu_1}^{2}-2\,{\it a_0}
\,b\mu_1-4\,{\it a_0}\,\mu_0-{b}^{2}+4\,{\it a_0} \right) }{8 \left( c
\sigma-{\it a_0} \right) }} \,,
\\[5 pt]
\bar \lambda= \sqrt{{\frac {3\,{{\it a_0}}^{2}{\mu_1}^{2}-2\,{\it a_0}\,b\mu_1-4
\,{\it a_0}\,\mu_0-{b}^{2}+4\,{\it a_0}}{16{\it a_0}\, \left( c\sigma-{
\it a_0} \right) }}}
 \,,
\end{cases} \displaybreak[3]\\
&\begin{cases}
\omega=- \left( {\it a_0}\,{\mu_1}^{2}-\mu_1\,b-\mu_0+\sigma \right) \mu_1 \,,
\\[5 pt]
\sigma > \frac{a_0}{2\,c} \,,
\end{cases}
\end{align*}
and one arrives at the solitary-wave solution 
\begin{equation*}
u(x,t) = e^{i\omega t} e^{iB(x-\sigma t)}[\bar d_0+ \bar f(x-\sigma t) ]
\ \ \ \ \ \mbox{and} \ \ \ \ \ v(x,t) =\bar h_0  \pm \sqrt{ \dd\frac{a_0}{2c\sigma-a_0}}\bar f(x-\sigma t)
\end{equation*} where $\bar f(\xi) = \bar d_2 \operatorname{sech}^2 (\bar \lambda \xi) $.

Aside from the above non-trivial solutions, system \eqref{system3} also possessess the following trivial and semi-trivial solutions
\begin{enumerate}
\item
\begin{equation*}
u(x,t) =0  \ \ \ \ \ \mbox{and} \ \ \ \ \ v(x,t) =h_0;
\end{equation*} for any $h_0 \in \mathbb R$.
\item
\begin{equation*}
u(x,t) = e^{i\omega t} e^{iB(x-\sigma t)} d_0  \ \ \ \ \ \mbox{and} \ \ \ \ \ v(x,t) =h_0;
\end{equation*} where $\s = \frac{\omega - a_0B^3 - b B^2 + B h_0 + B \mu_0 + h_0 \mu_1}{B}$, for any $B, d_0, h_0, \omega \in \mathbb R$.
\item
\begin{equation*}
u(x,t) = e^{i\omega t} e^{iB(x-\sigma t)}d_1  \cn\big(\lambda (x-\sigma t),m\big) \ \ \ \ \ \mbox{and} \ \ \ \ \ v(x,t) =h_0 + h_2  \cn^2\big(\lambda (x-\sigma t),m\big);
\end{equation*} 

for any $m \in [\,0,1\,]$, where $ {\it h_0}={\frac {9\,{{\it a_0}}^{2}c{m}^{2}{\mu_1}^{2}-6\,{\it a_0}\,bc{m}^{2}
      \mu_1-12\,{\it a_0}\,c{\it h_2}\,{m}^{2}-12\,{\it a_0}\,c{m}^{2}\mu_0-3\,{b}^{2}c{m}^{2}+2\,{
      {\it a_0}}^{2}{m}^{2}+6\,{\it a_0}\,c{\it h_2}}{12\,{\it a_0}\,c{m}^{2}}}$; \\
$h_2>0$ such that $ 9{{\it a_0}}^{2}{m}^{2}{\mu_1}^{2}-6{
      \it a_0}b{m}^{2}\mu_1
       -4{\it a_0}{\it h_2}{m}^{2}-12{\it a_0}{m}^{2}\mu_0  -3{b}^{2
      }{m}^{2}+2{\it a_0}{\it h_2}+12{\it a_0}{m}^{2}<0; $\\
$ d_1 = \pm {\frac {\sqrt {-6\,{\it a_0}h_2{m}^{2}\, \left( 9\,{{\it a_0}}^{2}
      {m}^{2}{\mu_1}^{2}-6\,{\it a_0}\,b{m}^{2}\mu_1-4\,{\it a_0}\,h_2\,{m}^{2}-12\,{\it a_0}
      \,{m}^{2}\mu_0-3\,{b}^{2}{m}^{2}+2\,{\it a_0}\,{\it h_2}+12\,{\it a_0}\,{m}^{2} \right) }}{6\,{
      \it a_0}\,{m}^{2}}}$;\\
$ B={\frac {{\it a_0}\,\mu_1-b}{2\,{\it a_0}}}$, $\omega =  {\frac {-\mu_1\, \left( 6\,{\it a_0}\,c{\mu_1}^{2}-6\,bc\mu_1-
      6\,c\mu_0+{\it a_0} \right) }{6\,c}}, $ $\l = \sqrt{\frac {{\it h_2}}{2\,{\it a_0}{m}^{2}}}$, and $\sigma={\frac {{\it a_0}}{6\,c}}$.
\item 
\begin{equation*}
u(x,t) =0\ \ \ \ \ \mbox{and} \ \ \ \ \ v(x,t) = -\frac{2}{3}\,{\it h_2}+\frac{1}{3}\,{\frac {{\it h_2}}{{m}^{2}}}+\sigma-1 + h_2  \cn^2\big(\lambda (x-\sigma t),m\big);
\end{equation*} where $ \l = \sqrt{\frac{h_2}{12c{m}^{2}\sigma}},$ for any $h_2, \sigma>0$ and $m \in [\,0,1\,]$.
\end{enumerate}

\subsection{Schr\"odinger BBM-KdV}
Setting all $k_{j,q}=0$ gives us the following set of parameters with $R$ as defined in \eqref{R}:
\begin{align*}
&\begin{cases}
B = {\frac {{\it a_1}\,\mu_0\,\mu_1-b}{2{\it a_1}\s\, \left( {\it a_1
}\,{\mu_1}^{2}+1 \right) }} \,,
\\[5 pt]
d_1 = h_1 = 0 \,,
\\[5 pt]
{\it d_0}= \frac{\sqrt{2\,c - {\it a_1}\,\sigma} \left( {m}^{4}-2\,{m}^{2}R-{m}^{2}+R+
1 \right) }{8\,R^2\sqrt{a_1\,\sigma}\left( {\it a_1}\,{\mu_1}^{2}+1 \right) ^{2} \left( {\it a_1
}\,\sigma-c \right)}\Big( 4\,{{\it a_1}}^{3}{\mu_1}^{4}\sigma-4\,{{
\it a_1}}^{2}b{\mu_1}^{3}\sigma-{{\it a_1}}^{2}{\mu_0}^{2}{\mu_1}^{2}-4\,{{
\it a_1}}^{2}\mu_0\,{\mu_1}^{2}\sigma +8\,{{\it a_1}}^{2}{\mu_1}^{2}\sigma
\\
\quad\quad +2
\,{\it a_1}\,b\mu_0\,\mu_1-4\,a_1b\mu_1\,\sigma-4\,a_1\mu_0\,\sigma+4\,{\it a_1}\,\sigma-{b}^{2} \Big)
 \,,
\\[5 pt]
{\it d_2}= \frac{3\,{m}^{2}\sqrt{2\,c - {\it a_1}\,\sigma}}{8R\sqrt{a_1\,\sigma}\left( {\it a_1}\,{\mu_1}^{2}+1 \right) ^{2} \left( {\it a_1
}\,\sigma-c \right) } \Big( 4\,{{\it a_1}}^{3}{\mu_1}^{4}\sigma-4\,{{
\it a_1}}^{2}b{\mu_1}^{3}\sigma-{{\it a_1}}^{2}{\mu_0}^{2}{\mu_1}^{2}-4\,{{
\it a_1}}^{2}\mu_0\,{\mu_1}^{2}\sigma+8\,{{\it a_1}}^{2}{\mu_1}^{2}\sigma
\\
\quad\quad +2
\,{\it a_1}\,b\mu_0\,\mu_1-4\,a_1b\mu_1\,\sigma-4\,a_1\mu_0\,\sigma
+4\,{\it a_1}\,\sigma-{b}^{2} \Big) \,,
 \\[5 pt]
{\it h_0}=\frac{-1}{8\,{\it a_1}\,R\sigma\, \left( {\it a_1}\,
{\mu_1}^{2}+1 \right) ^{2} \left( {\it a_1}\,\sigma - c \right)} \Big(-8\,{{\it a_1}}^{4}{\mu_1}^{4}R{\sigma}^{3}+8\,{{
\it a_1}}^{4}{m}^{2}{\mu_1}^{4}{\sigma}^{2}+4\,{{\it a_1}}^{4}{\mu_1}^{4}R{
\sigma}^{2}+8\,{{\it a_1}}^{3}c{\mu_1}^{4}R{\sigma}^{2}
\end{cases} \displaybreak[3]\\
&\begin{cases}
\quad\quad -4\,{{\it a_1}}^{4
}{\mu_1}^{4}{\sigma}^{2}-8\,{{\it a_1}}^{3}b{m}^{2}{\mu_1}^{3}{\sigma}^{2}+4\,{
{\it a_1}}^{3}b{\mu_1}^{3}R{\sigma}^{2}+4\,{{\it a_1}}^{3}b{\mu_1}^{3}{
\sigma}^{2}-2\,{{\it a_1}}^{3}{m}^{2}{\mu_0}^{2}{\mu_1}^{2}\sigma-8\,{{\it a_1}}
^{3}{m}^{2}\mu_0\,{\mu_1}^{2}{\sigma}^{2}
\\
\quad\quad +{{\it a_1}}^{3}{\mu_0}^{2}{\mu_1}^{2}
\sigma\,R+4\,{{\it a_1}}^{3}\mu_0\,{\mu_1}^{2}R{\sigma}^{2}-16\,{{\it a_1}
}^{3}{\mu_1}^{2}R{\sigma}^{3}-8\,a_1^2bc{\mu_1}^{3}R\sigma+16
\,{{\it a_1}}^{3}{m}^{2}{\mu_1}^{2}{\sigma}^{2}+{{\it a_1}}^{3}{\mu_0}^{2}{\mu_1}
^{2}\sigma
\\
\quad\quad +4\,{{\it a_1}}^{3}\mu_0\,{\mu_1}^{2}{\sigma}^{2}+8\,{{\it a_1}}
^{3}{\mu_1}^{2}R{\sigma}^{2}-2\,{{\it a_1}}^{2}c{\mu_0}^{2}{\mu_1}^{2}R-8\,
{{\it a_1}}^{2}c\mu_0\,{\mu_1}^{2}R\sigma+16\,{{\it a_1}}^{2}c{\mu_1}^{
2}R{\sigma}^{2}-8\,{{\it a_1}}^{3}{\mu_1}^{2}{\sigma}^{2}
\\
\quad\quad +4\,{{\it a_1}}^{2}b{m}^{2}\mu_0\,\mu
1\,\sigma-8\,{{\it a_1}}^{2}b{m}^{2}\mu_1\,{\sigma}^{2}-2\,{{\it a_1}}^{2}b
\mu_0\,\mu_1R\sigma+4\,{{\it a_1}}^{2}b\mu_1\,R{\sigma}^
{2}-2\,{{\it a_1}}^{2}b\mu_0\,\mu_1\,\sigma+4\,{{\it a_1}}^{2}b\mu_1\,{
\sigma}^{2}
\\
\quad\quad -8\,{{\it a_1}}^{2}{m}^{2}\mu_0\,{\sigma}^{2}+4\,{{\it a_1}}^{2}\mu_0
\,R{\sigma}^{2}-8\,{{\it a_1}}^{2}R{\sigma}^{3}+4\,{\it a_1}\,bc\mu_0\,
\mu_1\,R-8\,a_1 bc\mu_1\,R\sigma+8\,{{\it a_1}}^{2}{m}^{2}{\sigma}^{2}+4
\,{{\it a_1}}^{2}\mu_0\,{\sigma}^{2}
\\
\quad\quad +4\,{{\it a_1}}^{2}R{\sigma}^{2}-2\,a_1{
b}^{2}{m}^{2}\sigma+a_1{b}^{2}R\sigma-8\,a_1 c\mu_0\,R\sigma
+8\,{\it a_1}\,cR{\sigma}^{2}-4\,{{\it a_1}}^{2}{\sigma}^{2}+a_1{b}^
{2}\sigma-2\,{b}^{2}cR \Big) \,,
\\[5 pt]
{\it h_2}={\frac { 3\,\left( 4\,{{\it a_1}}^{3}{\mu_1}^{4}\sigma-4\,{{
\it a_1}}^{2}b{\mu_1}^{3}\sigma-{{\it a_1}}^{2}{\mu_0}^{2}{\mu_1}^{2}-4\,{{
\it a_1}}^{2}\mu_0\,{\mu_1}^{2}\sigma+8\,{{\it a_1}}^{2}{\mu_1}^{2}\sigma+2
\,{\it a_1}\,b\mu_0\,\mu_1-4\,a_1b\mu_1\,\sigma-4\,a_1\mu_0\,\sigma+4\,{\it 
a_1}\,\sigma-{b}^{2} \right) {m}^{2}}{8R\left( {\it a_1}\,{\mu_1
}^{2}+1 \right) ^{2} \left( {\it a_1}\,\sigma-c \right)}} \,,
\\[5 pt]
\lambda= \sqrt{\frac {4\,{{\it a_1}}^{3}{\mu_1}^{4}\sigma-4\,{{\it a_1}}^
{2}b{\mu_1}^{3}\sigma-{{\it a_1}}^{2}{\mu_0}^{2}{\mu_1}^{2}-4\,{{\it a_1}}^
{2}\mu_0\,{\mu_1}^{2}\sigma+8\,{{\it a_1}}^{2}{\mu_1}^{2}\sigma+2\,{\it a_1
}\,b\mu_0\,\mu_1-4\,a_1b\mu_1\,\sigma-4\,a_1\mu_0\,\sigma+4
\,{\it a_1}\,\sigma-{b}^{2}}{16\,{\it a_1}\,R\sigma\, \left( {\it a_1}\,{\mu_1}
^{2}+1 \right) ^{2} \left( {\it a_1}\,\sigma-c \right)}} \,,
\\[5 pt]
\omega=-{\frac { \left( {\it a_1}\,{\mu_1}^{2}\sigma-b\mu_1-\mu_0+\sigma
 \right) \mu_1}{{\it a_1}\,{\mu_1}^{2}+1}} \,,
\\[5 pt]
\sigma < \frac{2\,c}{a_1} \,,
\\[5 pt]
m \in [\,0,1\,] \,.
\end{cases}
\end{align*}
Thus, explicit periodic traveling-wave solutions to the Schr\"odinger BBM-KdV system $\big(u(x,t), v(x,t)\big) = \big(e^{i\omega t} e^{iB(x-\sigma t)}f(x-\sigma t), g(x-\sigma t)\big)$  given in term of the Jacobi cnoidal function

\begin{equation*}
f(\x)= d_0+d_2 \cn^2(\lambda \xi,m) \quad \text{and} \quad g(\x) = h_0 + h_2 \cn^2(\lambda \xi,m)
\end{equation*} are established.  Notice that $\frac{h_2}{d_2}= \sqrt{\frac{a_1\sigma}{2c-a_1\sigma}}$.  When $m=R=1$, the above coefficients simplify to $d_0=0$ and
$$
\begin{cases}
{\it \tilde d_2}= \frac{3\sqrt{2\,c - {\it a_1}\,\sigma}}{8\,\sqrt{a_1\,\sigma}\left( {\it a_1}\,{\mu_1}^{2}+1 \right) ^{2} \left( {\it a_1
}\,\sigma-c \right) } \Big( 4\,{{\it a_1}}^{3}{\mu_1}^{4}\sigma-4\,{{
\it a_1}}^{2}b{\mu_1}^{3}\sigma-{{\it a_1}}^{2}{\mu_0}^{2}{\mu_1}^{2}-4\,{{
\it a_1}}^{2}\mu_0\,{\mu_1}^{2}\sigma+8\,{{\it a_1}}^{2}{\mu_1}^{2}\sigma
\\
\quad\quad +2
\,{\it a_1}\,b\mu_0\,\mu_1-4\,a_1b\mu_1\,\sigma-4\,a_1\mu_0\,\sigma
+4\,{\it a_1}\,\sigma-{b}^{2} \Big) \,,
 \\[5 pt]
{\it \tilde h_0}=\frac{-1}{8\,{\it a_1}\,\sigma\, \left( {\it a_1}\,
{\mu_1}^{2}+1 \right) ^{2} \left( {\it a_1}\,\sigma - c \right) } \Big(-8\,{{\it a_1}}^{4}{\mu_1}^{4}{\sigma}^{3}+8\,{{
\it a_1}}^{4}{\mu_1}^{4}{\sigma}^{2}+8\,{{\it a_1}}^{3}c{\mu_1}^{4}{\sigma}^{2}
-16\,{{\it a_1}
}^{3}{\mu_1}^{2}{\sigma}^{3}-8{{\it a_1}}^{2}\,bc{\mu_1}^{3}\sigma\\
\quad\quad +16
\,{{\it a_1}}^{3}{\mu_1}^{2}{\sigma}^{2}
 -2\,{{\it a_1}}^{2}c{\mu_0}^{2}{\mu_1}^{2}-8
\,{{\it a_1}}^{2}c\mu_0\,{\mu_1}^{2}\sigma+16\,{{\it a_1}}^{2}c{\mu_1}^{
2}{\sigma}^{2}
-8\,{{\it a_1}}^{2}{\sigma}^{3}+4\,{\it a_1}\,bc\mu_0\,
\mu_1\\
\quad\quad-8{\it a_1}\,bc\mu_1\,\sigma+8\,{{\it a_1}}^{2}{\sigma}^{2}
-8{\it a_1}\,c\mu_0\,\sigma+8\,{\it a_1}\,c{\sigma}^{2}-2\,{b}^{2}c \Big) \,,
\\[5 pt]
{\it \tilde h_2}={\frac { 3\,\left( 4\,{{\it a_1}}^{3}{\mu_1}^{4}\sigma-4\,{{
\it a_1}}^{2}b{\mu_1}^{3}\sigma-{{\it a_1}}^{2}{\mu_0}^{2}{\mu_1}^{2}-4\,{{
\it a_1}}^{2}\mu_0\,{\mu_1}^{2}\sigma+8\,{{\it a_1}}^{2}{\mu_1}^{2}\sigma+2
\,{\it a_1}\,b\mu_0\,\mu_1-4\,a_1b\mu_1\,\sigma-4\,a_1\mu_0\,\sigma
+4\,{\it a_1}\,\sigma-{b}^{2} \right) }{8\, \left( {\it a_1}\,{\mu_1
}^{2}+1 \right) ^{2} \left( {\it a_1}\,\sigma-c \right) }} \,,
\\[5 pt]
\tilde \lambda= \sqrt{\frac {4\,{{\it a_1}}^{3}{\mu_1}^{4}\sigma-4\,{{\it a_1}}^
{2}b{\mu_1}^{3}\sigma-{{\it a_1}}^{2}{\mu_0}^{2}{\mu_1}^{2}-4\,{{\it a_1}}^
{2}\mu_0\,{\mu_1}^{2}\sigma+8\,{{\it a_1}}^{2}{\mu_1}^{2}\sigma+2\,{\it a_1
}\,b\mu_0\,\mu_1-4\,a_1b\mu_1\,\sigma-4\,a_1\mu_0\,\sigma+4
\,{\it a_1}\,\sigma-{b}^{2}}{16\,{\it a_1}\,\sigma\, \left( {\it a_1}\,{\mu_1}
^{2}+1 \right) ^{2} \left( {\it a_1}\,\sigma-c \right) }} \,,
\\[5 pt]
\omega=-{\frac { \left( {\it a_1}\,{\mu_1}^{2}\sigma-b\mu_1-\mu_0+\sigma
 \right) \mu_1}{{\it a_1}\,{\mu_1}^{2}+1}} \,,
\\[5 pt]
\sigma < \frac{2\,c}{a_1} \,,
\end{cases}
$$
from which one obtains the following solitary-wave solution to the system \eqref{system4} 
\begin{equation*}
u(x,t) = e^{i\omega t} e^{iB(x-\sigma t)} \tilde f(x-\sigma t) \ \ \ \ \ \mbox{and} \ \ \ \ \ v(x,t) =\tilde h_0  \pm \sqrt{\dd\frac{a_1\sigma}{2c-a_1 \sigma}} \tilde f(x-\sigma t)
\end{equation*}
where $\tilde f(\xi) = \tilde d_2 \operatorname{sech}^2 (\tilde \lambda \xi) $.  Furthermore, when $B = {\frac {{\it a_1}\,\mu_0\,\mu_1-b}{2\,\sigma\,{\it a_1}\, \left( {\it a_1
}\,{\mu_1}^{2}+1 \right) }}  $ satisfies the condition
\begin{align*}
    &(2 a_1^2 b c \mu_1^2 + 2 a_1 b c - 2 a_1^3 c \mu_0 \mu_1^3 - 2 a_1^2 c \mu_0 \mu_1) B^3 
    + (- 4 a_1^2 b c \mu_1^3 - 4 a_1^2 c \mu_0 \mu_1^2 - 4 a_1 b c \mu_1 - 4 a_1 c \mu_0) B^2 \notag \\
    &+ (2 a_1^3 \mu_0 \mu_1^3 + 2 a_1^2 c \mu_0 \mu_1^3 + 2 a_1^2 \mu_0 \mu_1 + 2 a_1 c \mu_0 \mu_1 - 2 a_1^2 b \mu_1^2 - 2 a_1 b - 2 a_1 b c \mu_1^2 - 2 b c) B \notag \\
    &+ (2 a_1 b \mu_0 \mu_1 - a_1^2 \mu_0^2 \mu_1^2 - b^2) = 0,
\end{align*}
one has $\tilde h_0=0$ and the synchronized solitary-wave solution established in \cite{NLB} is recovered.  

When $m=-R=1$, the above coefficients simplify to
$$
\begin{cases}
{\it \bar d_0}= \frac{\sqrt{2\,c - {\it a_1}\,\sigma}  }{4\,\sqrt{a_1\,\sigma}\left( {\it a_1}\,{\mu_1}^{2}+1 \right) ^{2} \left( {\it a_1
}\,\sigma-c \right) }\Big( 4\,{{\it a_1}}^{3}{\mu_1}^{4}\sigma-4\,{{
\it a_1}}^{2}b{\mu_1}^{3}\sigma-{{\it a_1}}^{2}{\mu_0}^{2}{\mu_1}^{2}-4\,{{
\it a_1}}^{2}\mu_0\,{\mu_1}^{2}\sigma +8\,{{\it a_1}}^{2}{\mu_1}^{2}\sigma
\\
\quad\quad +2
\,{\it a_1}\,b\mu_0\,\mu_1-4\,a_1b\mu_1\,\sigma-4\,a_1\mu_0\,\sigma
+4\,{\it a_1}\,\sigma-{b}^{2} \Big)
 \,,
\\[5 pt]
{\it \bar d_2}= \frac{-3\sqrt{2\,c - {\it a_1}\,\sigma}}{8\,\sqrt{a_1\,\sigma}\left( {\it a_1}\,{\mu_1}^{2}+1 \right) ^{2} \left( {\it a_1
}\,\sigma-c \right) } \Big( 4\,{{\it a_1}}^{3}{\mu_1}^{4}\sigma-4\,{{
\it a_1}}^{2}b{\mu_1}^{3}\sigma-{{\it a_1}}^{2}{\mu_0}^{2}{\mu_1}^{2}-4\,{{
\it a_1}}^{2}\mu_0\,{\mu_1}^{2}\sigma+8\,{{\it a_1}}^{2}{\mu_1}^{2}\sigma
\\
\quad\quad +2
\,{\it a_1}\,b\mu_0\,\mu_1-4\,a_1b\mu_1\,\sigma-4\,a_1\mu_0\,\sigma
+4\,{\it a_1}\,\sigma-{b}^{2} \Big) \,,
 \\[5 pt]
{\it \bar h_0}=\frac{1}{4\,{\it a_1}\,\sigma\, \left( {\it a_1}\,
{\mu_1}^{2}+1 \right) ^{2} \left( {\it a_1}\,\sigma - c \right) } \Big(4\,{{\it a_1}}^{4}{\mu_1}^{4}{\sigma}^{3}-4\,{{\it a_1}}^{3}c{\mu_1}^{4}{\sigma}^{2}
 -4\,{{\it a_1}}^{3}b{\mu_1}^{3}{\sigma}^{2}-{{\it a_1}}^{3}{\mu_0}^{2}{\mu_1}^{2}\sigma
\\
\quad\quad -4\,{{\it a_1}}
^{3}\mu_0\,{\mu_1}^{2}{\sigma}^{2} +8\,{{\it a_1}
}^{3}{\mu_1}^{2}{\sigma}^{3}+4{{\it a_1}}^{2}\,bc{\mu_1}^{3}\sigma
+{{\it a_1}}^{2}c{\mu_0}^{2}{\mu_1}^{2}+4{{\it a_1}}^{2}
\,c\mu_0\,{\mu_1}^{2}\sigma
\\
\quad\quad-8\,{{\it a_1}}^{2}c{\mu_1}^{
2}{\sigma}^{2}
+2{{\it a_1}}^{2}\,b\mu_0\,\mu_1
\,\sigma\,-4\,{{\it a_1}}^{2}b\mu_1\,{\sigma}^{2}
 -4\,{{\it a_1}}^{2}\mu_0\,{\sigma}^{2}+4\,{{\it a_1}}^{2}{\sigma}^{3}-2\,{\it a_1}\,bc\mu_0\,
\mu_1
\\
\quad\quad+4{\it a_1}\,bc\mu_1\,\sigma -{\it a_1}{b}^{2}\sigma +4{
\it a_1}\,c\mu_0\,\sigma-4\,{\it a_1}\,c{\sigma}^{2}+{b}^{2}c \Big) \,,
\\[5 pt]
{\it \bar h_2}={\frac { -3\,\left( 4\,{{\it a_1}}^{3}{\mu_1}^{4}\sigma-4\,{{
\it a_1}}^{2}b{\mu_1}^{3}\sigma-{{\it a_1}}^{2}{\mu_0}^{2}{\mu_1}^{2}-4\,{{
\it a_1}}^{2}\mu_0\,{\mu_1}^{2}\sigma+8\,{{\it a_1}}^{2}{\mu_1}^{2}\sigma+2
\,{\it a_1}\,b\mu_0\,\mu_1-4\,a_1b\mu_1\,\sigma-4\,a_1\mu_0\,\sigma
+4\,{\it a_1}\,\sigma-{b}^{2} \right) }{8\, \left( {\it a_1}\,{\mu_1
}^{2}+1 \right) ^{2} \left( {\it a_1}\,\sigma-c \right) }} \,,
\\[5 pt]
\bar \lambda= \sqrt{\frac {4\,{{\it a_1}}^{3}{\mu_1}^{4}\sigma-4\,{{\it a_1}}^
{2}b{\mu_1}^{3}\sigma-{{\it a_1}}^{2}{\mu_0}^{2}{\mu_1}^{2}-4\,{{\it a_1}}^
{2}\mu_0\,{\mu_1}^{2}\sigma+8\,{{\it a_1}}^{2}{\mu_1}^{2}\sigma+2\,{\it a_1
}\,b\mu_0\,\mu_1-4\,a_1b\mu_1\,\sigma-4\,a_1\mu_0\,\sigma+4
\,{\it a_1}\,\sigma-{b}^{2}}{-16\,{\it a_1}\,\sigma\, \left( {\it a_1}\,{\mu_1}
^{2}+1 \right) ^{2} \left( {\it a_1}\,\sigma-c \right) }} \,,
\\[5 pt]
\omega=-{\frac { \left( {\it a_1}\,{\mu_1}^{2}\sigma-b\mu_1-\mu_0+\sigma
 \right) \mu_1}{{\it a_1}\,{\mu_1}^{2}+1}} \,,
\\[5 pt]
\sigma < \frac{2\,c}{a_1} \,,
\end{cases}
$$
and one arrives at the solitary-wave solution 
\begin{equation*}
u(x,t) = e^{i\omega t} e^{iB(x-\sigma t)}[\bar d_0+ \bar f(x-\sigma t) ]
\ \ \ \ \ \mbox{and} \ \ \ \ \ v(x,t) =\bar h_0  \pm \sqrt{ \dd\frac{a_1\sigma}{2c-a_1\sigma}}\bar f(x-\sigma t)
\end{equation*} where $\bar f(\xi) = \bar d_2 \operatorname{sech}^2 (\bar \lambda \xi) $.

Aside from the above non-trivial solutions, system \eqref{system4} also possessess the following trivial and semi-trivial solutions 
\begin{enumerate}
\item
\begin{equation*}
u(x,t) =0  \ \ \ \ \ \mbox{and} \ \ \ \ \ v(x,t) =h_0;
\end{equation*} for any $h_0 \in \mathbb R$.
\item
\begin{equation*}
u(x,t) = e^{i\omega t} e^{iB(x-\sigma t)} d_0  \ \ \ \ \ \mbox{and} \ \ \ \ \ v(x,t) = h_0;
\end{equation*} where $\s = \frac{a_1 B^2 \omega - b B^2 + B h_0 + B \mu_0 + h_0 \mu_1 + \omega}{B (a_1 B^2 + 1)} $, for any $B,d_0,h_0,\omega \in \mathbb R$.
\item
\begin{equation*}
u(x,t) = e^{i\omega t} e^{iB(x-\sigma t)} d_1 \cn\big(\lambda (x-\sigma t),m\big) \ \ \ \ \ \mbox{and} \ \ \ \ \ v(x,t) = h_0+ h_2 \cn^2\big(\lambda (x-\sigma t),m\big);
\end{equation*} where $d_1 = \pm \frac{1}{6\,c{m}^{2}(a_1 \mu_1^2 + 1)}
      \Big( 3\,c{\it h_2}{m}^{2}\, \big( 8\,{{\it a_1}}^{2}c{\it h_2}\,{m}^{2}{\mu_1}^{4}-4
      \,{{\it a_1}}^{2}c{\it h_2}{\mu_1}^{4}-24\,{{\it a_1}}^{2}c{m}^{2}{\mu_1}^{4}+{{\it a_1}}
      ^{2}{m}^{2}{\mu_0}^{2}{\mu_1}^{2} +24\,{\it a_1}\,bc{m}^{2}{\mu_1}^{3}
       +16\,{\it a_1}\,c{
      \it h_2}\,{m}^{2}{\mu_1}^{2}+24\,{\it a_1}\,c{m}^{2}\mu_0\,{\mu_1}^{2}-2\,{\it a_1}\,b{m}^{2}
      \mu_0\,\mu_1-8\,{\it a_1}\,c{\it h_2}{\mu_1}^{2}-48\,{\it a_1}\,c{m}^{2}{\mu_1}^{
      2}+24\,bc{m}^{2}\mu_1+{b}^{2}{m}^{2}
       +8\,c{\it h_2}\,{m}^{2}+24\,c{m}^{2}\mu_0-4\,c{\it h_2}-24\,
      c{m}^{2} \big) \Big)^{1/2}$;\\
$ h_0 = \frac{-1}{24\, \left( {{\it a_1}}^{2}{\mu_1}^{4}+2
      \,a_1 {\mu_1}^{2}+1 \right) {\it a_1}\,c{m}^{2}}
      \Big( 24\,{{\it a_1}}^{3}c{\it h_2}\,{m}^{2}{\mu_1}^{4}-12\,{
      {\it a_1}}^{3}c{\it h_2}\,{\mu_1}^{4}-144\,{{\it a_1}}^{2}{c}^{2}{m}^{2}{\mu_1}^{
      4}+{{\it a_1}}^{3}{m}^{2}{\mu_0}^{2}{\mu_1}^{2}
      \\\quad\quad +24\,{{\it a_1}}^{2}bc{m}^{2}{\mu_1}^{3}
      +48\,{{\it a_1}}^{2}c{\it h_2}\,{m}^{2}{\mu_1}^{2}+24\,{{\it a_1}}^{2}c{m}^{2}\mu_0\,{
      \mu_1}^{2}-2\,{{\it a_1}}^{2}b{m}^{2}\mu_0\,\mu_1-24\,{{\it a_1}}^{2}c{\it h_2}\,{
      \mu_1}^{2}
      \\\quad\quad -288\,{\it a_1}\,{c}^{2}{m}^{2}{\mu_1}^{2}+24\,{\it a_1}\,bc{m}^{2}\mu_1+{
      \it a_1}\,{b}^{2}{m}^{2}+24\,{\it a_1}\,c{\it h_2}\,{m}^{2}+24\,{\it a_1}\,c{m}^{2}\mu_0-12\,
      {\it a_1}\,c{\it h_2}-144\,{c}^{2}{m}^{2} \Big)$;\\
any $h_2>0$ such that
      $8\,{{\it a_1}}^{2}c{\it h_2}\,{m}^{2}{\mu_1}^{4}-4
      \,{{\it a_1}}^{2}c{\it h_2}{\mu_1}^{4}-24\,{{\it a_1}}^{2}c{m}^{2}{\mu_1}^{4}+{{\it a_1}}
      ^{2}{m}^{2}{\mu_0}^{2}{\mu_1}^{2} +24\,{\it a_1}\,bc{m}^{2}{\mu_1}^{3}
      +16\,{\it a_1}\,c{
      \it h_2}\,{m}^{2}{\mu_1}^{2}+24\,{\it a_1}\,c{m}^{2}\mu_0\,{\mu_1}^{2}-2\,{\it a_1}\,b{m}^{2}
      \mu_0\,\mu_1-8\,{\it a_1}\,c{\it h_2}{\mu_1}^{2}-48\,{\it a_1}\,c{m}^{2}{\mu_1}^{
      2}+24\,bc{m}^{2}\mu_1+{b}^{2}{m}^{2}
      +8\,c{\it h_2}\,{m}^{2}+24\,c{m}^{2}\mu_0-4\,c{\it h_2}-24\,
      c{m}^{2} >0$;\\
$ B= {\frac {{\it a_1}\,\mu_0\,\mu_1-b}{12\,c \left( a_1 {\mu_1}^{2}+1
      \right) }};$ $ \l = \sqrt{\frac{h_2}{12\,c{m}^{2}}};$ $\omega= {\frac {\mu_1\, \left( -6\,{\it a_1}\,c{\mu_1}^{2}+a_1 b \mu_1
      +{\it a_1}\,\mu_0-6\,c \right) }{ \left( a_1 {\mu_1}^{2}+1 \right) {
      \it a_1}}};$ $\sigma= {\frac {6\,c}{{\it a_1}}};$
 and any $ m \in [\,0,1\,]$.
\item
\begin{equation*}
u(x,t) = 0 \ \ \ \ \ \mbox{and} \ \ \ \ \ v(x,t) = -\frac{2}{3}\,{\it h_2}+\frac{1}{3}\,{\frac {{\it h_2}}{{m}^{2}}}+\sigma-1 + h_2\cn^2\big(\lambda (x-\sigma t),m\big);
\end{equation*} where $  \l = \sqrt{\frac{h_2}{12c{m}^{2}}} $, for any $h_2,\sigma>0$ and $m \in [\,0,1\,]$.
\end{enumerate}

\section{Conclusion} The periodic traveling wave solutions for the four systems \eqref{system1}-\eqref{system4} are found.  Our results show that \textit{all} periodic solutions to the four systems are given by \eqref{travelingwavesolution} and \eqref{cnoidalsolutions}.  These cnoidal solutions limit to the solitary-wave solutions when $m \rightarrow 1$.  This is expected since it is well-known that the ODE equation 
\begin{equation*}
f'^2 = k_3f^3 +k_2f^2 + k_1 f + k_0
\end{equation*}
has unique solitary-wave solution as well as periodic cnoidal solution, and that the periodic cnoidal solution limits to the solitary-wave solution when the Jacobi elliptic modulus $m$ approaches $1$.  Consequently, we obtain solitary-wave solutions for all four systems as the by-products.  All of the \textit{synchronized} solitary-wave solutions established in \cite{NLB} are special cases of those obtained here.   Another direct consequence is that the \textit{synchronized} periodic solutions previously obtained in \cite{BLN} are indeed unique, a fact that wasn't established therein.  Since those synchronized periodic solutions approach the synchronized solitary-wave solutions obtained in \cite{NLB}, it would be interesting to know whether these synchronized solitary-wave solutions are also indeed unique.  This question is not pursued here.  

Below are some graphs for the cnoidal wave solutions 
for the four systems \eqref{system1}-\eqref{system4}.  
Recall that a traveling-wave solution to the above four 
systems is a vector solution $\big(u(x,t),v(x,t)\big)$ 
of the form 
\begin{equation*}
u(x,t) = e^{i\omega t} e^{iB(x-\sigma t)}f(x-\sigma t),\ \ \ \ \ \ \ v(x,t) = g(x-\sigma t),
\end{equation*}
where $f$ and $g$ are smooth, real-valued functions 
with speed $\sigma>0$, and phase shifts $B,\omega 
\in \mathbb R$. For ease of graphing, the imaginary terms in $u(x,t)$
were supressed as they define a phase shift and thus a rotation 
of the real function $f$, which is graphed below. For all four vector solutions,
$m=\frac{1}{2}$ and $R = \frac{\sqrt{13}}{4}$ were chosen, while the 
remaining parameters were then fixed to ensure real solutions
and are listed here; KdV-KdV: 
$\sigma=2, a_0=1, b=-1, c=\frac{3}{2}, \mu_0=1, \mu_1=\frac{1}{4}$; BBM-BBM: $\sigma=1, 
a_1=1, b=-1, c=\frac{5}{2}, \mu_0=1, \mu_1=1$; KdV-BBM:
$\sigma=\frac{3}{2}, a_0=1, b=-1, c=\frac{3}{2}, \mu_0=1, \mu_1=\frac{1}{4}$;
BBM-KdV: $\sigma=\frac{1}{2}, a_1=1, b=-1, c=\frac{3}{2}, \mu_0=1, \mu_1=\frac{1}{4}$.
The graphs are now listed below, with $u(x,t)$ in blue and 
$v(x,t)$ in green.
\begin{figure}[h]
  \centering
  \begin{subfigure}[b]{0.4\linewidth}
    \includegraphics[width=\linewidth]{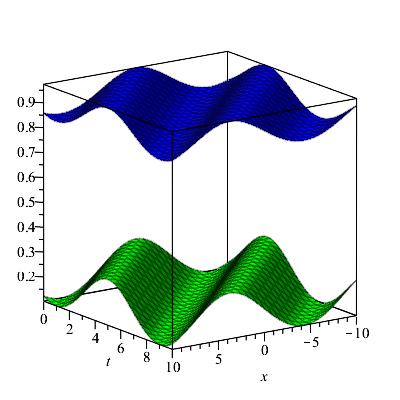}
    \caption{Cnoidal solution for the KdV-KdV system.}
  \end{subfigure}
  \begin{subfigure}[b]{0.4\linewidth}
    \includegraphics[width=\linewidth]{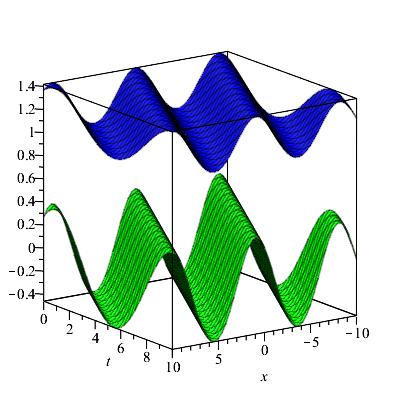}
    \caption{Cnoidal solution for the BBM-BBM system.}
  \end{subfigure}
\begin{subfigure}[b]{0.4\linewidth}
    \includegraphics[width=\linewidth]{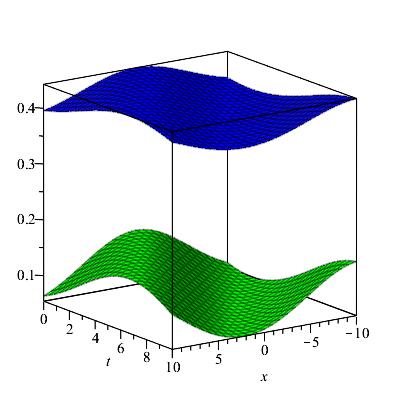}
    \caption{Cnoidal solution for the KdV-BBM system.}
  \end{subfigure}
  \begin{subfigure}[b]{0.4\linewidth}
    \includegraphics[width=\linewidth]{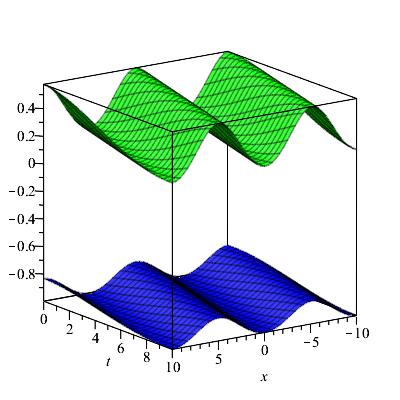}
    \caption{Cnoidal solution for the BBM-KdV system.}
  \end{subfigure}
\end{figure}

\section{Appendix}
For the Schr\"odinger KdV-KdV system \eqref{system1}, the $k_{j, q}$ in \eqref{genericsystem} are:
\begin{equation} \label{coeffKdVKdV}
\begin{cases}
k_{1,3} = 4\,{\lambda}^{2}\,{\it a_0}\,{\it d_2}\,{m}^{2} - \frac{2}{3}\,{\it d_2}\,{\it h_2} \,,
\\[5 pt]
k_{1,2} = {\lambda}^{2}\,{\it a_0}\,{\it d_1}\,{m}^{2} - \frac{1}{2}\,{\it d_1}\,{\it h_2} - \frac{1}{2}\,{\it d_2}\,
{\it h_1} \,,
\\[5 pt]
k_{1,1} = {\it d_2}\,{\it a_0}\,{B}^{2}-\frac{8}{3}\,{\lambda}^{2}\,{\it a_0}\,{\it d_2}\,{m}^{2}+\frac{2}{3}\,{
\it d_2}\,bB+\frac{4}{3}\,{\it d_2}\,{\it a_0}\,{\lambda}^{2}-\frac{1}{3}\,{\it d_0}\,{\it h_2}-\frac{1}{3}\,{\it d_1}\,{\it h_1}
\\ 
\quad\quad -\frac{1}{3}\,{\it d_2}\,{\it h_0}-\frac{1}{3}\,{\it d_2}\,\mu_0+\frac{1}{3}\,{\it d_2}\,\sigma \,,
\end{cases}$$
$$\begin{cases}
k_{1,0} = \frac{1}{2}\,{B}^{2}{\it a_0}\,{\it d_1}-\frac{1}{3}\,{\lambda}^{2}\,{\it a_0}\,{\it d_1}\,{m}^{2}+\frac{1}{3}
\,Bb{\it d_1}+\frac{1}{6}\,{\lambda}^{2}\,{\it a_0}\,{\it d_1}-\frac{1}{6}\,{\it d_0}\,{\it h_1}
\\
\quad\quad -\frac{1}{6}\,{\it d_1}\,{\it h_0}-\frac{1}{6}\,{\it d_1}\,\mu_0+\frac{1}{6}\,{\it d_1}\,\sigma \,,
\\[5 pt]
k_{2,4} = -18\,B{\it a_0}\,{\it d_2}\,{\lambda}^{2}\,{m}^{2}-6\,b{\it d_2}\,{\lambda}^{2}\,{m}^{2}+B{\it d_2
}\,{\it h_2}+{\it d_2}\,{\it h_2}\,\mu_1 \,, 
\\[5 pt]
k_{2,3} = -6\,B{\it a_0}\,{\it d_1}\,{\lambda}^{2}\,{m}^{2}-2\,b{\it d_1}\,{\lambda}^{2}\,{m}^{2}+B{\it d_1}
\,{\it h_2}+B{\it d_2}\,{\it h_1}+{\it d_1}\,{\it h_2}\,\mu_1+{\it d_2}\,{
\it h_1}\,\mu_1 \,, \\[5 pt]
k_{2,2} = -{B}^{3}{\it a_0}\,{\it d_2}+24\,B{\it a_0}\,{\it d_2}\,{\lambda}^{2}\,{m}^{2}-{B}^{2}
b{\it d_2}-12\,B{\it a_0}\,{\it d_2}\,{\lambda}^{2}+8\,b{\it d_2}\,{\lambda}^{2}\,{m}^{2}+B{
\it d_0}\,{\it h_2}
\\
\quad\quad +B{\it d_1}\,{\it h_1}+B{\it d_2}\,{\it h_0}+B{\it d_2}\,
\mu_0-B{\it d_2}\,\sigma-4\,b{\it d_2}\,{\lambda}^{2}+{\it d_0}\,{\it h_2}\,\mu_1+
{\it d_1}\,{\it h_1}\,\mu_1
\\
\quad\quad +{\it d_2}\,{\it h_0}\,\mu_1+{\it d_2}\,\omega \,,
\\[5 pt]
k_{2,1} = -{B}^{3}{\it a_0}\,{\it d_1}+6\,B{\it a_0}\,{\it d_1}\,{\lambda}^{2}\,{m}^{2}-{B}^{2}b
{\it d_1}-3\,B{\it a_0}\,{\it d_1}\,{\lambda}^{2}+2\,b{\it d_1}\,{\lambda}^{2}\,{m}^{2}+B{
\it d_0}\,{\it h_1}
\\
\quad\quad +B{\it d_1}\,{\it h_0}+B{\it d_1}\,\mu_0-B{\it d_1}\,
\sigma-b{\it d_1}\,{\lambda}^{2}+{\it d_0}\,{\it h_1}\,\mu_1+{\it d_1}\,{\it h_0}
\,\mu_1+{\it d_1}\,\omega \,,
\\[5 pt]
k_{2,0} = -{B}^{3}{\it a_0}\,{\it d_0}-6\,B{\it a_0}\,{\it d_2}\,{\lambda}^{2}\,{m}^{2}-{B}^{2}b
{\it d_0}+6\,B{\it a_0}\,{\it d_2}\,{\lambda}^{2}-2\,b{\it d_2}\,{\lambda}^{2}\,{m}^{2}+B{
\it d_0}\,{\it h_0}
\\
\quad\quad +B{\it d_0}\,\mu_0-B{\it d_0}\,\sigma+2\,b{\it d_2}\,
{\lambda}^{2}+{\it d_0}\,{\it h_0}\,\mu_1+{\it d_0}\,\omega \,,
\\[5 pt]
k_{3,3} = \frac{1}{12}\,{{\it d_2}}^{2}+\frac{1}{12}\,{{\it h_2}}^{2}-{\lambda}^{2}\,c{\it h_2}\,{m}^{2} \,,
\\[5 pt]
k_{3,2} = -\frac{1}{4}\,{\lambda}^{2}\,c{\it h_1}\,{m}^{2}+\frac{1}{8}\,{\it d_1}\,{\it d_2}+\frac{1}{8}\,{\it h_1}\,{
\it h_2} \,,
\\[5 pt]
k_{3,1} = -\frac{1}{3}\,{\lambda}^{2}\,c{\it h_2}+\frac{2}{3}\,{\lambda}^{2}\,c{\it h_2}\,{m}^{2}+\frac{1}{12}\,{\it d_0}\,{
\it d_2}+\frac{1}{12}\,{\it h_0}\,{\it h_2}-\frac{1}{12}\,{\it h_2}\,\sigma+\frac{1}{24}\,{{\it d_1
}}^{2}
\\
\quad\quad +\frac{1}{24}\,{{\it h_1}}^{2}+\frac{1}{12}{\it h_2} \,,
\\[5 pt]
k_{3,0} = \frac{1}{12}\,{\lambda}^{2}\,c{\it h_1}\,{m}^{2}+\frac{1}{24}\,{\it d_0}\,{\it d_1}+\frac{1}{24}\,{\it h_0}\,{
\it h_1}-\frac{1}{24}\,{\it h_1}\,\sigma-\frac{1}{24}\,{\lambda}^{2}\,c{\it h_1}+\frac{1}{24}{\it h_1} \,.
\end{cases}
\end{equation}

For the Schr\"odinger BBM-BBM system \eqref{system2}, the $k_{j, q}$ in \eqref{genericsystem} are:
\begin{equation} \label{coeffBBMBBM}
\begin{cases}
k_{1,3} = -\frac{2}{3}\,{\it d_2}\,{\it h_2}+4\,{\lambda}^{2}\,{\it a_1}\,{\it d_2}\,{m}^{2}\sigma \,,
\\[5 pt]
k_{1,2} = -\frac{1}{2}\,{\it d_1}\,{\it h_2}-\frac{1}{2}\,{\it d_2}\,{\it h_1}+{\lambda}^{2}\,{\it a_1}\,{
\it d_1}\,{m}^{2}\sigma \,,
\\[5 pt]
k_{1,1} = \frac{2}{3}\,Bb{\it d_2}-\frac{1}{3}\,{\it d_0}\,{\it h_2}-\frac{1}{3}\,{\it d_1}\,{\it h_1}-\frac{1}{3}\,{
\it d_2}\,{\it h_0}-\frac{1}{3}\,{\it d_2}\,\mu_0+\frac{1}{3}\,{\it d_2}\,\sigma
\\
\quad\quad +\frac{4}{3}\,
{\lambda}^{2}\,{\it a_1}\,{\it d_2}\,\sigma+{B}^{2}{\it a_1}\,{\it d_2}\,\sigma-
\frac{2}{3}\,B{\it a_1}\,{\it d_2}\,\omega-\frac{8}{3}\,{\lambda}^{2}\,{\it a_1}\,{\it d_2}\,{m}^{2}
\sigma \,,
\\[5 pt]
k_{1,0} = \frac{1}{3}\,Bb{\it d_1}-\frac{1}{6}\,{\it d_0}\,{\it h_1}-\frac{1}{6}\,{\it d_1}\,{\it h_0}-\frac{1}{6}\,{
\it d_1}\,\mu_0+\frac{1}{6}\,{\it d_1}\,\sigma+\frac{1}{2}\,{B}^{2}{\it a_1}\,{\it d_1}\,
\sigma
\\
\quad\quad -\frac{1}{3}\,B{\it a_1}\,{\it d_1}\,\omega-\frac{1}{3}\,{\lambda}^{2}\,{\it a_1}\,{\it 
d_1}\,{m}^{2}\sigma+\frac{1}{6}\,{\lambda}^{2}\,{\it a_1}\,{\it d_1}\,\sigma \,,
\\[5 pt]
k_{2,4} = -18\,B{\it a_1}\,{\it d_2}\,{\lambda}^{2}\,{m}^{2}\sigma+6\,{\it a_1}\,{\it d_2}\,
{\lambda}^{2}\,{m}^{2}\omega-6\,b{\it d_2}\,{\lambda}^{2}\,{m}^{2}+B{\it d_2}\,{\it h_2}+{\it d_2}
\,{\it h_2}\,\mu_1 \,,
\\[5 pt]
k_{2,3} = -6\,B{\it a_1}\,{\it d_1}\,{\lambda}^{2}\,{m}^{2}\sigma+2\,{\it a_1}\,{\it d_1}\,
{\lambda}^{2}\,{m}^{2}\omega-2\,b{\it d_1}\,{\lambda}^{2}\,{m}^{2}+B{\it d_1}\,{\it h_2}+B{\it d_2
}\,{\it h_1}
\\
\quad\quad +{\it d_1}\,{\it h_2}\,\mu_1+{\it d_2}\,{\it h_1}\,\mu_1 \,,
\\[5 pt]
k_{2,2} = 6\,B{\it a_1}\,{\it d_2}\,{\lambda}^{2}\,{m}^{2}\sigma+6\, \left( -B \left( -{m}^{2}+1
 \right) +B{m}^{2} \right) {\it a_1}\,{\it d_2}\,{\lambda}^{2}\,\sigma-6\,B \left( -
{m}^{2}+1 \right) {\it a_1}\,{\it d_2}\,{\lambda}^{2}\,\sigma
\\
\quad\quad -2\,{\it a_1}\,{\it d_2}
\,{\lambda}^{2}\,{m}^{2}\omega-{B}^{3}{\it a_1}\,{\it d_2}\,\sigma+{B}^{2}{\it a_1}\,
{\it d_2}\,\omega-2\, \left( 2\,{m}^{2}-1 \right) {\it a_1}\,{\it d_2}\,{\lambda}^{2}
\,\omega
\\
\quad\quad +2\, \left( -{m}^{2}+1 \right) {\it a_1}\,{\it d_2}\,{\lambda}^{2}\,\omega+2
\,b{\it d_2}\,{\lambda}^{2}\,{m}^{2}+B{\it d_2}\,\mu_0+B{\it d_0}\,{\it h_2}+{\it d_1}\,
{\it h_1}\,\mu_1+{\it d_2}\,{\it h_0}\,\mu_1
\\
\quad\quad -{B}^{2}b{\it d_2}+{\it d_2}\,
\omega-B{\it d_2}\,\sigma+{\it d_0}\,{\it h_2}\,\mu_1+B{\it d_1}\,{\it h_1}+
B{\it d_2}\,{\it h_0}+2\, \left( 2\,{m}^{2}-1 \right) b{\it d_2}\,{\lambda}^{2}
\\
\quad\quad -2\,
 \left( -{m}^{2}+1 \right) b{\it d_2}\,{\lambda}^{2} \,,
\end{cases}$$
$$\begin{cases}
k_{2,1} = 3\,B{\it a_1}\,{\it d_1}\,{\lambda}^{2}\,{m}^{2}\sigma-3\,B \left( -{m}^{2}+1 \right) {
\it a_1}\,{\it d_1}\,{\lambda}^{2}\,\sigma-{\it a_1}\,{\it d_1}\,{\lambda}^{2}\,{m}^{2}
\omega+{B}^{2}{\it a_1}\,{\it d_1}\,\omega
\\
\quad\quad -{B}^{3}{\it a_1}\,{\it d_1}\,
\sigma+ \left( -{m}^{2}+1 \right) {\it a_1}\,{\it d_1}\,{\lambda}^{2}\,\omega+b{\it 
d_1}\,{\lambda}^{2}\,{m}^{2}-{B}^{2}b{\it d_1}+B{\it d_0}\,{\it h_1}
\\
\quad\quad +B{\it d_1}\,{\it 
h_0}+B{\it d_1}\,\mu_0-B{\it d_1}\,\sigma+{\it d_0}\,{\it h_1}\,\mu_1+{\it d_1
}\,{\it h_0}\,\mu_1+{\it d_1}\,\omega- \left( -{m}^{2}+1 \right) b{\it d_1}\,
{\lambda}^{2} \,,
\\[5 pt]
k_{2,0} = -{B}^{3}{\it a_1}\,{\it d_0}\,\sigma+{B}^{2}{\it a_1}\,{\it d_0}\,\omega+{
\it d_0}\,\omega-{B}^{2}b{\it d_0}+B{\it d_0}\,{\it h_0}+B{\it d_0}\,\mu_0-B
{\it d_0}\,\sigma
\\
\quad\quad +{\it d_0}\,{\it h_0}\,\mu_1+6\,B \left( -{m}^{2}+1 \right) {
\it a_1}\,{\it d_2}\,{\lambda}^{2}\,\sigma-2\, \left( -{m}^{2}+1 \right) {\it a_1}\,{
\it d_2}\,{\lambda}^{2}\,\omega+2\, \left( -{m}^{2}+1 \right) b{\it d_2}\,{\lambda}^{2} \,,
\\[5 pt]
k_{3,3} = 24\,c{\it h_2}\,{\lambda}^{2}\,{m}^{2}\sigma-2\,{{\it d_2}}^{2}-2\,{{\it h_2}}^{2} \,,
\\[5 pt]
k_{3,2} = 6\,{\lambda}^{2}\,c{\it h_1}\,{m}^{2}\sigma-3\,{\it d_1}\,{\it d_2}-3\,{\it h_1}\,{
\it h_2} \,, \\[5 pt]
k_{3,1} = -16\,c{\it h_2}\,{\lambda}^{2}\,{m}^{2}\sigma+8\,c{\it h_2}\,{\lambda}^{2}\,\sigma-2\,{
\it d_0}\,{\it d_2}-{{\it d_1}}^{2}-2\,{\it h_0}\,{\it h_2}-{{\it h_1}}^{2}+
2\,{\it h_2}\,\sigma-2\,{\it h_2} \,,
\\[5 pt]
k_{3,0} = -2\,{\lambda}^{2}\,c{\it h_1}\,{m}^{2}\sigma+{\lambda}^{2}\,c{\it h_1}\,\sigma-{\it d_0}\,{
\it d_1}-{\it h_0}\,{\it h_1}+{\it h_1}\,\sigma-{\it h_1} \,.
\end{cases}
\end{equation}

For the Schr\"odinger KdV-BBM system \eqref{system3}, the $k_{j, q}$ in \eqref{genericsystem} are:
\begin{equation} \label{coeffKdVBBM}
      \begin{cases}
          k_{1,3} = 4\,{\lambda}^{2}\,{\it a_0}\,{\it d_2}\,{m}^{2} - \frac{2}{3}\,{\it d_2}\,{\it h_2} \,,
          \\[5 pt]
          k_{1,2} = {\lambda}^{2}\,{\it a_0}\,{\it d_1}\,{m}^{2} - \frac{1}{2}\,{\it d_1}\,{\it h_2} - \frac{1}{2}\,{\it d_2}\,
      {\it h_1} \,,
      \\[5 pt]
      k_{1,1} = {\it d_2}\,{\it a_0}\,{B}^{2}-\frac{8}{3}\,{\lambda}^{2}\,{\it a_0}\,{\it d_2}\,{m}^{2}+\frac{2}{3}\,{
      \it d_2}\,bB+\frac{4}{3}\,{\it d_2}\,{\it a_0}\,{\lambda}^{2}-\frac{1}{3}\,{\it d_0}\,{\it h_2}-\frac{1}{3}\,{\it d_1}\,{\it h_1}
      \\ 
      \quad\quad -\frac{1}{3}\,{\it d_2}\,{\it h_0}-\frac{1}{3}\,{\it d_2}\,\mu_0+\frac{1}{3}\,{\it d_2}\,\sigma \,,
      \\[5 pt]
      k_{1,0} = \frac{1}{2}\,{B}^{2}{\it a_0}\,{\it d_1}-\frac{1}{3}\,{\lambda}^{2}\,{\it a_0}\,{\it d_1}\,{m}^{2}+\frac{1}{3}
      \,Bb{\it d_1}+\frac{1}{6}\,{\lambda}^{2}\,{\it a_0}\,{\it d_1}-\frac{1}{6}\,{\it d_0}\,{\it h_1}
      \\
      \quad\quad -\frac{1}{6}\,{\it d_1}\,{\it h_0}-\frac{1}{6}\,{\it d_1}\,\mu_0+\frac{1}{6}\,{\it d_1}\,\sigma \,,
      \\[5 pt]
      k_{2,4} = -18\,B{\it a_0}\,{\it d_2}\,{\lambda}^{2}\,{m}^{2}-6\,b{\it d_2}\,{\lambda}^{2}\,{m}^{2}+B{\it d_2
      }\,{\it h_2}+{\it d_2}\,{\it h_2}\,\mu_1 \,,
      \\[5 pt]
      k_{2,3} = -6\,B{\it a_0}\,{\it d_1}\,{\lambda}^{2}\,{m}^{2}-2\,b{\it d_1}\,{\lambda}^{2}\,{m}^{2}+B{\it d_1}
      \,{\it h_2}+B{\it d_2}\,{\it h_1}+{\it d_1}\,{\it h_2}\,\mu_1+{\it d_2}\,{
      \it h_1}\,\mu_1 \,,
      \\[5 pt]
      k_{2,2} = -{B}^{3}{\it a_0}\,{\it d_2}+24\,B{\it a_0}\,{\it d_2}\,{\lambda}^{2}\,{m}^{2}-{B}^{2}
      b{\it d_2}-12\,B{\it a_0}\,{\it d_2}\,{\lambda}^{2}+8\,b{\it d_2}\,{\lambda}^{2}\,{m}^{2}+B{
      \it d_0}\,{\it h_2}
      \\
      \quad\quad +B{\it d_1}\,{\it h_1}+B{\it d_2}\,{\it h_0}+B{\it d_2}\,
      \mu_0-B{\it d_2}\,\sigma-4\,b{\it d_2}\,{\lambda}^{2}+{\it d_0}\,{\it h_2}\,\mu_1+
      {\it d_1}\,{\it h_1}\,\mu_1
      \\
      \quad\quad +{\it d_2}\,{\it h_0}\,\mu_1+{\it d_2}\,\omega \,,
      \\[5 pt]
      k_{2,1} = -{B}^{3}{\it a_0}\,{\it d_1}+6\,B{\it a_0}\,{\it d_1}\,{\lambda}^{2}\,{m}^{2}-{B}^{2}b
      {\it d_1}-3\,B{\it a_0}\,{\it d_1}\,{\lambda}^{2}+2\,b{\it d_1}\,{\lambda}^{2}\,{m}^{2}+B{
      \it d_0}\,{\it h_1}
      \\
      \quad\quad +B{\it d_1}\,{\it h_0}+B{\it d_1}\,\mu_0-B{\it d_1}\,
      \sigma-b{\it d_1}\,{\lambda}^{2}+{\it d_0}\,{\it h_1}\,\mu_1+{\it d_1}\,{\it h_0}
      \,\mu_1+{\it d_1}\,\omega \,,
\\[5 pt]
k_{2,0} = -{B}^{3}{\it a_0}\,{\it d_0}-6\,B{\it a_0}\,{\it d_2}\,{\lambda}^{2}\,{m}^{2}-{B}^{2}b
{\it d_0}+6\,B{\it a_0}\,{\it d_2}\,{\lambda}^{2}-2\,b{\it d_2}\,{\lambda}^{2}\,{m}^{2}+B{
\it d_0}\,{\it h_0} 
\\
\quad\quad +B{\it d_0}\,\mu_0-B{\it d_0}\,\sigma+2\,b{\it d_2}\,
{\lambda}^{2}+{\it d_0}\,{\it h_0}\,\mu_1+{\it d_0}\,\omega \,,
\\[5 pt]
k_{3,3} = \frac{1}{3}\,{{\it d_2}}^{2}+\frac{1}{3}\,{{\it h_2}}^{2}-4\,{\lambda}^{2}\,c{\it h_2}\,{m}^{2}\sigma \,,
\\[5 pt]
k_{3,2} = \frac{1}{2}\,{\it d_1}\,{\it d_2}+\frac{1}{2}\,{\it h_1}\,{\it h_2}-{\lambda}^{2}\,c{\it h_1}\,{m}^{2}
\sigma \,,
\\[5 pt]
k_{3,1} = \frac{1}{3}\,{\it d_0}\,{\it d_2}+\frac{1}{3}\,{\it h_0}\,{\it h_2}-\frac{1}{3}\,{\it h_2}\,\sigma-
\frac{4}{3}\,{\lambda}^{2}\,c{\it h_2}\,\sigma+\frac{1}{6}\,{{\it d_1}}^{2}+\frac{1}{6}\,{{\it h_1}}^{2
}+\frac{1}{3}{\it h_2} + \frac{8}{3}\,{\lambda}^{2}\,c{\it h_2}\,{m}^{2}\sigma \,,
\\[5 pt]
k_{3,0} = \frac{1}{6}{\it h_1}-\frac{1}{6}\,{\lambda}^{2}\,c{\it h_1}\,\sigma+\frac{1}{6}\,{\it d_0}\,{\it d_1}+\frac{1}{6}
\,{\it h_0}\,{\it h_1}-\frac{1}{6}\,{\it h_1}\,\sigma+\frac{1}{3}\,{\lambda}^{2}\,c{\it h_1}\,{m}^{2}
\sigma \,.
\end{cases}\end{equation}
\newpage

For the Schr\"odinger BBM-KdV system \eqref{system4}, the $k_{j, q}$ in \eqref{genericsystem} are:
\begin{equation} \label{coeffBBMKdV}
\begin{cases}
k_{1,3} = -\frac{2}{3}\,{\it d_2}\,{\it h_2}+4\,{\lambda}^{2}\,{\it a_1}\,{\it d_2}\,{m}^{2}\sigma \,,
\\[5 pt]
k_{1,2} = -\frac{1}{2}\,{\it d_1}\,{\it h_2}-\frac{1}{2}\,{\it d_2}\,{\it h_1}+{\lambda}^{2}\,{\it a_1}\,{
\it d_1}\,{m}^{2}\sigma \,,
\\[5 pt]
k_{1,1} = \frac{2}{3}\,Bb{\it d_2}-\frac{1}{3}\,{\it d_0}\,{\it h_2}-\frac{1}{3}\,{\it d_1}\,{\it h_1}-\frac{1}{3}\,{
\it d_2}\,{\it h_0}-\frac{1}{3}\,{\it d_2}\,\mu_0+\frac{1}{3}\,{\it d_2}\,\sigma
\\
\quad\quad +\frac{4}{3}\,
{\lambda}^{2}\,{\it a_1}\,{\it d_2}\,\sigma+{B}^{2}{\it a_1}\,{\it d_2}\,\sigma-
\frac{2}{3}\,B{\it a_1}\,{\it d_2}\,\omega-\frac{8}{3}\,{\lambda}^{2}\,{\it a_1}\,{\it d_2}\,{m}^{2}
\sigma \,,
\\[5 pt]
k_{1,0} = \frac{1}{3}\,Bb{\it d_1}-\frac{1}{6}\,{\it d_0}\,{\it h_1}-\frac{1}{6}\,{\it d_1}\,{\it h_0}-\frac{1}{6}\,{
\it d_1}\,\mu_0+\frac{1}{6}\,{\it d_1}\,\sigma+\frac{1}{2}\,{B}^{2}{\it a_1}\,{\it d_1}\,
\sigma
\\
\quad\quad -\frac{1}{3}\,B{\it a_1}\,{\it d_1}\,\omega-\frac{1}{3}\,{\lambda}^{2}\,{\it a_1}\,{\it 
d_1}\,{m}^{2}\sigma+\frac{1}{6}\,{\lambda}^{2}\,{\it a_1}\,{\it d_1}\,\sigma \,,
\\[5 pt]
k_{2,4} = -18\,B{\it a_1}\,{\it d_2}\,{\lambda}^{2}\,{m}^{2}\sigma+6\,{\it a_1}\,{\it d_2}\,
{\lambda}^{2}\,{m}^{2}\omega-6\,b{\it d_2}\,{\lambda}^{2}\,{m}^{2}+B{\it d_2}\,{\it h_2}+{\it d_2}
\,{\it h_2}\,\mu_1 \,,
\\[5 pt]
k_{2,3} = -6\,B{\it a_1}\,{\it d_1}\,{\lambda}^{2}\,{m}^{2}\sigma+2\,{\it a_1}\,{\it d_1}\,
{\lambda}^{2}\,{m}^{2}\omega-2\,b{\it d_1}\,{\lambda}^{2}\,{m}^{2}+B{\it d_1}\,{\it h_2}+B{\it d_2
}\,{\it h_1}
\\
\quad\quad +{\it d_1}\,{\it h_2}\,\mu_1+{\it d_2}\,{\it h_1}\,\mu_1 \,,
\\[5 pt]
k_{2,2} = 6\,B{\it a_1}\,{\it d_2}\,{\lambda}^{2}\,{m}^{2}\sigma+6\, \left( -B \left( -{m}^{2}+1
 \right) +B{m}^{2} \right) {\it a_1}\,{\it d_2}\,{\lambda}^{2}\,\sigma-6\,B \left( -
{m}^{2}+1 \right) {\it a_1}\,{\it d_2}\,{\lambda}^{2}\,\sigma
\\
\quad\quad -2\,{\it a_1}\,{\it d_2}
\,{\lambda}^{2}\,{m}^{2}\omega-{B}^{3}{\it a_1}\,{\it d_2}\,\sigma+{B}^{2}{\it a_1}\,
{\it d_2}\,\omega-2\, \left( 2\,{m}^{2}-1 \right) {\it a_1}\,{\it d_2}\,{\lambda}^{2}
\,\omega
\\
\quad\quad +2\, \left( -{m}^{2}+1 \right) {\it a_1}\,{\it d_2}\,{\lambda}^{2}\,\omega+2
\,b{\it d_2}\,{\lambda}^{2}\,{m}^{2}+B{\it d_2}\,\mu_0+B{\it d_0}\,{\it h_2}+{\it d_1}\,
{\it h_1}\,\mu_1+{\it d_2}\,{\it h_0}\,\mu_1 
\\
\quad\quad -{B}^{2}b{\it d_2}+{\it d_2}\,
\omega-B{\it d_2}\,\sigma+{\it d_0}\,{\it h_2}\,\mu_1+B{\it d_1}\,{\it h_1}+
B{\it d_2}\,{\it h_0}+2\, \left( 2\,{m}^{2}-1 \right) b{\it d_2}\,{\lambda}^{2}
\\
\quad\quad -2\,
 \left( -{m}^{2}+1 \right) b{\it d_2}\,{\lambda}^{2} \,,
\\[5 pt]
k_{2,1} = 3\,B{\it a_1}\,{\it d_1}\,{\lambda}^{2}\,{m}^{2}\sigma-3\,B \left( -{m}^{2}+1 \right) {
\it a_1}\,{\it d_1}\,{\lambda}^{2}\,\sigma-{\it a_1}\,{\it d_1}\,{\lambda}^{2}\,{m}^{2}
\omega+{B}^{2}{\it a_1}\,{\it d_1}\,\omega
\\
\quad\quad -{B}^{3}{\it a_1}\,{\it d_1}\,
\sigma+ \left( -{m}^{2}+1 \right) {\it a_1}\,{\it d_1}\,{\lambda}^{2}\,\omega+b{\it 
d_1}\,{\lambda}^{2}\,{m}^{2}-{B}^{2}b{\it d_1}+B{\it d_0}\,{\it h_1}
\\
\quad\quad +B{\it d_1}\,{\it 
h_0}+B{\it d_1}\,\mu_0-B{\it d_1}\,\sigma+{\it d_0}\,{\it h_1}\,\mu_1+{\it d_1
}\,{\it h_0}\,\mu_1+{\it d_1}\,\omega- \left( -{m}^{2}+1 \right) b{\it d_1}\,
{\lambda}^{2} \,,
\\[5 pt]
k_{2,0} = -{B}^{3}{\it a_1}\,{\it d_0}\,\sigma+{B}^{2}{\it a_1}\,{\it d_0}\,\omega+{
\it d_0}\,\omega-{B}^{2}b{\it d_0}+B{\it d_0}\,{\it h_0}+B{\it d_0}\,\mu_0-B
{\it d_0}\,\sigma
\\
\quad\quad +{\it d_0}\,{\it h_0}\,\mu_1+6\,B \left( -{m}^{2}+1 \right) {
\it a_1}\,{\it d_2}\,{\lambda}^{2}\,\sigma-2\, \left( -{m}^{2}+1 \right) {\it a_1}\,{
\it d_2}\,{\lambda}^{2}\,\omega+2\, \left( -{m}^{2}+1 \right) b{\it d_2}\,{\lambda}^{2} \,,
\\[5 pt]
k_{3,3} = \frac{1}{12}\,{{\it d_2}}^{2}+\frac{1}{12}\,{{\it h_2}}^{2}-{\lambda}^{2}\,c{\it h_2}\,{m}^{2} \,,
\\[5 pt]
k_{3,2} = -\frac{1}{4}\,{\lambda}^{2}\,c{\it h_1}\,{m}^{2}+\frac{1}{8}\,{\it d_1}\,{\it d_2}+\frac{1}{8}\,{\it h_1}\,{
\it h_2} \,,
\\[5 pt]
k_{3,1} = -\frac{1}{3}\,{\lambda}^{2}\,c{\it h_2}+\frac{2}{3}\,{\lambda}^{2}\,c{\it h_2}\,{m}^{2}+\frac{1}{12}\,{\it d_0}\,{
\it d_2}+\frac{1}{12}\,{\it h_0}\,{\it h_2}-\frac{1}{12}\,{\it h_2}\,\sigma+\frac{1}{24}\,{{\it d_1
}}^{2} 
\\
\quad\quad +\frac{1}{24}\,{{\it h_1}}^{2}+\frac{1}{12}{\it h_2} \,,
\\[5 pt]
k_{3,0} = \frac{1}{12}\,{\lambda}^{2}\,c{\it h_1}\,{m}^{2}+\frac{1}{24}\,{\it d_0}\,{\it d_1}+\frac{1}{24}\,{\it h_0}\,{
\it h_1}-\frac{1}{24}\,{\it h_1}\,\sigma-\frac{1}{24}\,{\lambda}^{2}\,c{\it h_1}+\frac{1}{24}{\it h_1} \,.
\end{cases}
\end{equation}

\section*{Use of AI tools declaration}
The authors declare they have not used Artificial Intelligence (AI) tools in the creation of this article.
\section*{Conflict of interest} All authors declare no conflicts of interest in this paper.

\end{document}